\documentclass{article}

\usepackage{amssymb}
\usepackage{amsmath}
\usepackage[all]{xy}
\usepackage[frenchb]{babel}
\usepackage{a4wide}
\usepackage{enumerate,fancyhdr}

\begin{document}

\title{Intersection de courbes et de sous-groupes, et probl\`emes de minoration de hauteur dans les vari\'et\'es ab\'eliennes C.M.}

\author{Nicolas Ratazzi \footnote{nicolas.ratazzi@math.u-psud.fr\ \ \ \textbf{Adresse : } Universit\'e Paris-Sud 11, Batiment 425 Math\'ematiques, 91405 Orsay Cedex, France}}

\renewcommand{\text}{\textnormal}
\renewcommand{\hat}{\widehat}
\newcommand{\h}{\hat{h}}

\renewcommand{\phi}{\varphi}
\renewcommand{\epsilon}{\varepsilon}
\renewcommand{\tilde}{\widetilde}

\newcommand{\Gal}{\textnormal{Gal}}
\newcommand{\tors}{\textnormal{tors}}
\newcommand{\ord}{\textnormal{ord}}

\newcommand{\Z}{\mathbb{Z}}
\newcommand{\Q}{\mathbb{Q}}
\newcommand{\C}{\mathbb{C}}
\newcommand{\R}{\mathbb{R}}
\newcommand{\N}{\mathbb{N}}
\newcommand{\G}{\mathbb{G}}

\renewcommand{\S}{\Sigma}

\renewcommand{\O}{\mathcal{O}}
\renewcommand{\L}{\mathcal{L}}
\renewcommand{\P}{\mathcal{P}}
\newcommand{\p}{\mathfrak{p}}
\newcommand{\pp}{\mathfrak{P}}
\newcommand{\Dif}{\mathcal{D}}

\newcommand{\mmid}{\|}

\newcommand{\Deg}{\textnormal{Deg\,}}
\newcommand{\e}{\varepsilon}
\newcommand{\ab}{\textnormal{ab}}
\newcommand{\dd}{\delta}
\newcommand{\ddd}{D_{\text{inc}}}
\newcommand{\dds}{\delta^{\star}}
\newcommand{\Id}{\textnormal{Id}}
\renewcommand{\b}{\beta}
\newcommand{\bb}{\tilde{\beta}}
\newcommand{\tF}{\tilde{F}}
\renewcommand{\a}{\alpha}
\renewcommand{\aa}{\tilde{\alpha}}

\setcounter{tocdepth}{2}

\newcounter{ndefinition}[section]
\newcommand{\defi}{\addtocounter{ndefinition}{1}{\noindent \textbf{D{\'e}finition \thesection.\thendefinition.\ }}}
\newcounter{nrem}
\newcommand{\rem}{\addtocounter{nrem}{1}{\noindent \textbf{Remarque \thesection.\thenrem.\ }}}
\newtheorem{lemme}{Lemme}[section]
\newtheorem{conj}{Conjecture}[section]
\newcounter{nex}[section]
\newcommand{\ex}{\addtocounter{nex}{1}{\noindent \textit{Exemple} \thesection.\thenex\ }}
\newtheorem{theme}{Th{\`e}me} [section]
\newtheorem{prop}{Proposition} [section]
\newtheorem{cor}{Corollaire} [section]
\newtheorem{theo}{Th{\'e}or{\`e}me} [section]
\newcommand{\demo}{\noindent \textit{D{\'e}monstration} : }

\newcounter{compteur}

\makeatletter
\renewcommand{\thecompteur}{\@Alph\c@compteur}
\makeatother

\newtheorem{propo}[compteur]{Proposition}
\newtheorem{theor}[compteur]{Th\'eor\`eme}
\newtheorem{lem}[compteur]{Lemme}

\maketitle

\noindent \textbf{Abstract }We prove a special case of the following conjecture of Zilber-Pink generalising the Manin-Mumford conjecture : let $X$ be a curve inside an Abelian variety $A$ over $\overline{\Q}$, provided $X$ is not contained in a torsion subvariety, the intersection of $X$ with the union of all subgroup schemes of codimension at least $2$ is finite ; we settle the case where $A$ is a power of a simple Abelian variety of C.M. type. This generalises the previous known result, due to Viada and R\'emond-Viada (who was able to prove the conjecture for power of an elliptic curve with complex multiplication). The proof is based on the strategy of R\'emond (following Bombieri, Masser and Zannier) with two new ingredients, one of them, being at the heart of this article : it is a lower bound for the N\'eron-Tate height of points on Abelian varieties $A/K$ of C.M. type in the spirit of Lehmer's problem. This lower bound is an analog of the similar result of Amoroso and David \cite{ad2003} on $\G_m^n$ and is a generalisation of the theorem of David and Hindry \cite{davidhindry} on the abelian Lehmer's problem. The proof is an adaptation of \cite{davidhindry} using in our abelian case the new ideas introduced in \cite{ad2003}. Furthermore, as in \cite{ad2003} and adapting in the abelian case their proof, we give another application of our result : a lower bound for the absolute minimum of a subvariety $V$ of $A$. Although lower bounds for this minimum were already known (decreasing  multi-exponential function of the degree for Bombieri-Zannier), our methods enable us to prove, up to an $\varepsilon$ the optimal result that can be conjectured.

\bigskip

\tableofcontents

\bigskip

\noindent \textit{classification :} 11G50, 11G10, 11J95, 14K22, 11R20

\medskip

\noindent \textit{keywords : }Abelian varieties, normalised height, Lehmer's problem, Manin-Mumford Conjecture

\section{Introduction}

\subsection{Introduction}

\noindent Dans cet article nous nous int\'eressons \`a la g\'en\'eralisation suivante de la conjecture de Manin-Mumford : soit $G$ une vari\'et\'e semi-ab\'elienne sur $\overline{\Q}$ et $X$ une courbe contenue dans $G$ ; soit $G^{[r]}$ l'union de tous les sous-groupes alg\'ebriques non n\'ecessairement connexes de $G$ de codimension $\geq r$, on demande pour quelles valeurs de $r$ l'intersection $X\cap G^{[r]}$ est finie. La conjecture de Manin-Mumford correspond au cas $r=\dim G$ (en effet $G^{[\dim G]}=G(\overline{\Q})_{\tors})$. De plus pour avoir un $r$ aussi petit que possible, il faut visiblement supposer que $X$ n'est contenue dans aucun sous-groupe alg\'ebrique strict de $G$. L'\'enonc\'e est faux pour $r=1$ et la conjecture optimiste est $r=2$.

\medskip

\noindent Nous obtenons ici (\textit{cf.} th\'eor\`eme \ref{ray}), en suivant une strat\'egie d\'evelopp\'ee par Bombieri-Masser-Zannier \cite{BMZ}, Viada \cite{viada} puis R\'emond \cite{remond}, le r\'esultat optimal $r=2$ pour les vari\'et\'es ab\'eliennes de la forme $A=B^n$ o\`u $B$ est une vari\'et\'e ab\'elienne de type C.M. simple.  La preuve est notamment bas\'ee sur deux raffinements : d'une part une am\'elioration (faisant l'objet de la pr\'epublication s\'epar\'ee \cite{matorsion}) dans le cas des vari\'et\'es ab\'eliennes de type C.M., d'un r\'esultat de Masser \cite{lettre} concernant le probl\`eme de la borne (non-uniforme) sur la torsion ; d'autre part un raffinement du r\'esultat minoration de hauteur principal de David et Hindry \cite{davidhindry} sur le probl\`eme de Lehmer pour les vari\'et\'es ab\'eliennes de type C.M. C'est la preuve, par des techniques d'approximation diophantienne ou de transcendance, de ce dernier raffinement qui est au coeur du pr\'esent article.

\medskip

\noindent Ce type de probl\`eme (intersection de courbes et de sous groupes) a tout d'abord \'et\'e trait\'e dans \cite{BMZ} dans le cas de $\G_m^n$ pour lequel ils obtiennent un r\'esultat essentiellement optimal : r=2 mais leur hypoth\`ese sur $C$ est l\'eg\`erement plus forte. Leur m\'ethode (et l'obtention du r\'esultat $r=2$) a ensuite \'et\'e \'etendue par Viada \cite{viada}, compl\'et\'e par R\'emond-Viada \cite{rv}, au cas d'une vari\'et\'e ab\'elienne de la forme $A=E^n$ avec $E$ une courbe elliptique \`a multiplication complexe. R\'emond \cite{remond} a finalement \'etendu la strat\'egie (mais pas le r\'esultat optimal) au cas d'une vari\'et\'e ab\'elienne quelconque. N\'eanmoins il n'obtient pas de nouveau r\'esultat optimal inconditionnel (son r\'esultat incondionnel \'etant loin de l'optimalit\'e, voir par exemple le th\'eor\`eme \ref{gr} ci-dessous pour le cas CM). Par contre il montre qu'une tr\`es bonne minoration, conjecturale, de la hauteur des points non de torsion (\textit{cf.} la conjecture \ref{conj2} ci-apr\`es ainsi que la remarque qui suit) entra\^ine $r=2$ pour $G=A$ vari\'et\'e ab\'elienne.

\medskip

\noindent  Ainsi notre r\'esultat permet de passer du cas d'une puissance d'une courbe elliptique de type C.M. \`a une puissance d'une vari\'et\'e ab\'elienne de type C.M., simple de dimension quelconque. Avant de donner les \'enonc\'es pr\'ecis, indiquons d\`es \`a pr\'esent que notre \'enonc\'e cl\'e de minoration de hauteur (th\'eor\`eme \ref{semi}) est l'analogue dans le cas des vari\'et\'es ab\'eliennes de type C.M. du r\'esultat de Amoroso-David \cite{ad2003} dans le cas de $\G_m^n$. De plus notre preuve est une combinaison des preuves de \cite{davidhindry} et \cite{ad2003}~: nous reprenons le fil de la preuve de \cite{davidhindry} en introduisant dans notre cadre ab\'elien les id\'ees nouvelles de \cite{ad2003}.

\medskip

\noindent Par ailleurs, comme il est fait dans \cite{ad2003} pour $\G_m^n$ et suivant leur preuve nous tirons du r\'esultat de minoration de hauteur une seconde application (\textit{cf.} th\'eor\`eme \ref{absolu}), elle aussi li\'ee aux probl\`emes de minoration de hauteurs dans les vari\'et\'es ab\'eliennes. Finalement nous indiquons en appendice deux derniers r\'esultats~ concernant les probl\`emes de minoration de hauteur. Le premier donne une preuve du fait (qui est pr\'ecis\'e en appendice) que ''une bonne minoration de la hauteur des points d'ordre infini modulo toute sous-vari\'et\'e ab\'elienne sur une vari\'et\'e ab\'elienne entra\^ine une bonne minoration de la hauteur de tous les points d'ordre infini". Le second est une preuve de ce que le probl\`eme de Lehmer ab\'elien de \cite{davidhindry} (\textit{cf.} leur conjecture 1.4) est \'equivalent au probl\`eme \textit{a priori} plus fort, dit multihomog\`ene (\textit{cf.} la conjecture 1.6 de \cite{davidhindry}). Dans le cas de $\G_m^n$ cette \'equivalence avait d\'ej\`a \'et\'e not\'ee et d\'emontr\'ee par Amoroso et David dans leur article \cite{amodavlehm}. Nous adaptons en appendice leur preuve au cas des vari\'et\'es ab\'eliennes.

\subsection{\'Enonc\'es}

\subsubsection{Probl\`emes de minoration de hauteur}

\noindent Le r\'esultat que nous obtenons (th\'eor\`eme \ref{semi}) g\'en\'eralise le r\'esultat principal de \cite{davidhindry}. Il s'agit de l'analogue du m\^eme r\'esultat pour $\G_m^n$ d\'emontr\'e dans \cite{ad2003}. Notre preuve suit la preuve de \cite{davidhindry} en utilisant dans notre cadre ab\'elien les id\'ees nouvelles introduites dans \cite{ad2003}. Par ailleurs, tout comme cela \'etait fait dans \cite{amodavlehm} pour le cas de $\G_m^n$ nous rempla\c{c}ons dans le r\'esultat le degr\'e $[K(x):K]$ par \textit{l'indice d'obstruction} (\textit{cf.} d\'efinition \ref{delta}). Une telle possibilit\'e avait d\'ej\`a \'et\'e not\'ee dans \cite{davidhindry}. L'application la plus frappante de ce r\'esultat concerne les probl\`emes d'intersection de courbes et de sous-groupes alg\'ebriques et est d\'etaill\'ee plus bas dans la sous-section suivante. Enfin tout comme \cite{ad2003} dans le cas de $\G_m^n$ et en adaptant leur preuve au cadre ab\'elien nous donnons une autre application de notre r\'esultat : le th\'eor\`eme \ref{semi} concernant la minoration du minimum absolu des sous-vari\'et\'es de vari\'et\'es ab\'eliennes de type C.M. Dans ce qui suit nous utilisons librement la notion d'indice d'obstruction $\dd_{\L}$ rappel\'ee plus loin dans l'introduction.

\medskip

\begin{conj}\label{conj1app2}\textnormal{\textbf{(Probl\`eme de Lehmer ab\'elien)}} Soient $A/K$ une vari\'et\'e ab\'elienne sur un corps de nombres et $\L$ un fibr\'e en droites sym\'etrique ample sur $A$. Il existe une constante strictement positive $c(A/K,\L)$ telle que pour tout point $P\in A(\overline{K})$ d'ordre infini modulo toute sous-vari\'et\'e ab\'elienne stricte de $A$, on a
\begin{equation}\label{e1}
\widehat{h}_{\L}(P)\geq\frac{c(A/K,\L)}{\delta_{\L,K}(P)}.
\end{equation}
\noindent De plus, en terme du degr\'e $D=[K(P):K]$, on a pour tout point $P\in A(\overline{K})$ qui n'est pas de torsion
\begin{equation}\label{e2}
\widehat{h}_{\L}(P)\geq\frac{c(A/K,\L)}{D^{\frac{1}{g_0}}},
\end{equation}
\noindent o\`u $g_0$ est la dimension du plus petit sous-groupe alg\'ebrique contenant le point $P$.
\end{conj}

\noindent Dans cette direction,  David et Hindry obtiennent le r\'esultat suivant (c'est le th\'eor\`eme 1.5 de \cite{davidhindry}) :

\medskip

\begin{theor} \label{t1}\textnormal{\textbf{(David-Hindry \cite{davidhindry})}} Soient $A/K$ une vari\'et\'e ab\'elienne de type C.M. de dimension $g$ sur un corps de nombres et munie d'un fibr\'e en droites ample et sym\'etrique $\L$. Il existe une constante strictement positive $c(A/K,\L)$ telle que pour tout point $P\in A(\overline{K})$ d'ordre infini modulo toute sous-vari\'et\'e ab\'elienne stricte de $A$ on a 
\[ \h_{\L}(P)\geq \frac{c(A/K,\L)}{D^{\frac{1}{g}}}\left(\frac{\log \log 3D}{\log 2D}\right)^{\kappa(g)},\]
\noindent o\`u $\kappa(g)$ est une fonction explicite de $g$ et $D=[K(P):K]$.
\end{theor}

\medskip

\noindent Par ailleurs, ils \'enoncent un probl\`eme de Lehmer ``multihomog\`ene'' \textit{a priori} plus fort que le probl\`eme de Lehmer ab\'elien (\textit{cf.} le paragraphe \ref{multi} de l'appendice du pr\'esent article). Ils se demandent \'egalement dans quelle mesure leur th\'eor\`eme \ref{t1} pourrait \^etre raffin\'e afin de dire des choses sur la deuxi\`eme partie du probl\`eme de Lehmer ab\'elien. Enfin, ils indiquent qu'il serait int\'eressant de quantifier l'hypoth\`ese ``d'ordre infini'' en terme du degr\'e de la plus petite sous-vari\'et\'e de torsion pouvant contenir le point.

\medskip

\noindent Dans notre pr\'esent article nous r\'epondons entre autres \`a toutes ces questions, dans un cadre raffin\'e plus pr\'ecis. Nous formulons notre r\'esultat principal (et la preuve) en utilisant des indices d'obstructions. Avant d'\'enoncer pr\'ecis\'ement nos r\'esultats, commen\c{c}ons par d\'efinir l'indice d'obstruction suivant \cite{davidhindry} d\'efinition 1.2~:

\medskip

\defi \label{delta}Soient $A/K$ une vari\'et\'e ab\'elienne, $\L$ un fibr\'e en droites ample et sym\'etri\-que, $P$ un point de $A(\overline{K})$ et $F/K$ une extension alg\'ebrique. On d\'efinit \textit{l'indice d'obstruction de $P$ relativement \`a $\L$ et $F$}, et l'on note $\dd_{\L,F}(P)$, par 
\[\dd_{\L,F}(P)=\min\left\{\left(\deg_{\L_F}X\right)^{\frac{1}{\text{codim}X}} |\ X \text{ sous-$F$-vari\'et\'e stricte de $A_F$, telle que }P\in X(\overline{K})\right\},\]
\noindent o\`u l'on a not\'e $\L_F$ le faisceau sur $A_F$ tir\'e en arri\`ere de $\L$ par la projection naturelle de $A_F$ sur $A$. Plus g\'en\'eralement on peut d\'efinir l'\textit{indice d'obstruction (relativement \`a $\L$ et $F$), pour une sous}-$\overline{K}$-\textit{vari\'et\'e $V$ de }$A_{\overline{K}}$ :
\[\dd_{\L,F}(V)=\min\left\{\left(\deg_{\L_F}X\right)^{\frac{1}{\text{codim}X}}\ |\ \ X \ \text{sous-$F$-vari\'et\'e stricte de $A_F$, telle que }\overline{V}^F\subset X\right\},\]
\noindent o\`u l'on a not\'e $\overline{V}^F$ l'image sch\'ematique de $V\subset A_{\overline{K}}$ dans $A_F$.

\medskip

\rem \label{rem1}En consid\'erant la vari\'et\'e $\overline{\{P\}}^F$, image sch\'ematique de $P\in A_{\overline{K}}$ dans $A_F$, on constate imm\'ediatement que 
\[\dd_{\L,F}(P)\leq [F(P):F]^{\frac{1}{g}}.\]

\medskip

\noindent On notera $\dd_{\L,\text{tors}}$ l'indice d'obstruction relatif au corps $K_{\text{tors}}=K(A_{\tors})$. Par ailleurs dans la suite, comme on travaille avec un fibr\'e en droites $\L$ fixe, on l'omettra r\'egulierement dans la notation de $\dd_{\L,F}$ afin de ne pas trop alourdir les notations. Ceci \'etant, on peut maintenant \'enoncer le raffinement attendu du probl\`eme de Lehmer :

\medskip

\begin{conj} \label{conj2}\textnormal{\textbf{(Probl\`eme de Lehmer ab\'elien relatif)}} Soient $A/K$ une vari\'et\'e ab\'elienne sur un corps de nombres et $\L$ un fibr\'e en droites ample et sym\'etrique. Il existe une constante strictement positive $c(A/K,\L)$ telle que pour tout point $P\in A(\overline{K})$ d'ordre infini modulo toute sous-vari\'et\'e ab\'elienne stricte de $A$ on a 
\[ \h_{\L}(P)\geq \frac{c(A/K,\L)}{\dd_{\L,\tors}(P)}.\]
\end{conj}

\medskip

\rem Cette conjecture n'englobe \textit{a priori} que la premi\`ere partie de l'\'enonc\'e du probl\`eme de Lehmer ab\'elien de David et Hindry. Mais comme nous le montrons dans l'appendice (corollaire \ref{co2}), la seconde partie du probl\`eme de Lehmer ab\'elien de \cite{davidhindry} est en fait une cons\'equence de la premi\`ere.

\medskip

\noindent Dans un pr\'ec\'edent travail, \cite{art3}, nous avons obtenu dans le cas de ce probl\`eme de Lehmer relatif un r\'esultat, optimal \`a des puissances de $\log D_{\tors}$ pr\`es, dans le cas des courbes elliptiques \`a multiplication complexe. En direction de cette conjecture, nous obtenons, dans le cas des vari\'et\'es ab\'elien\-nes de type C.M., un r\'esultat, essentiellement optimal pour le terme principal. Pour tout entier $n$, on note $K_n:=K(A[n])$ l'extension engendr\'ee sur $K$ par le groupe des points de torsion $A[n]$. On a $K_{\tors}=\bigcup_{n\geq 1}K(A[n])$. 

\medskip

\begin{theo} \label{semi} Soit $A/K$ une vari\'et\'e ab\'elienne de type C.M. de dimension $g$ sur un corps de nombres et munie d'un fibr\'e en droites ample et sym\'etrique $\L$. Il existe une constante strictement positive $c(A/K,\L)$ telle que pour tout point $P\in A(\overline{K})$ et pour tout entier $n$, on a l'alternative suivante :
\[ \textit{soit }\hspace{1cm} \h_{\L}(P)\geq \frac{c(A/K,\L)}{\dd_{\L, K_n}(P)}\left(\frac{\log\log 3\left[K_n:K\right]\dd_{\L,K_n}(P)}{\log 2\left[K_n:K\right]\dd_{\L,K_n}(P)}\right)^{\kappa(g)},\]
\noindent avec  $\kappa(g)=(g+1)!(2g+5)(g+2)(2g.g!)^g$ ;

\noindent soit le point $P$ est contenu dans une sous-vari\'et\'e de torsion stricte, $B$, de $A_{K_n}$, d\'efinie sur $K_n$, de degr\'e major\'e par
\[\left(\deg_{\L_{K_n}} B\right)^{\frac{1}{\text{codim }B}}\leq c(A/K,\L)^{-1}\dd_{\L,K_n}(P)\left(\log 2\left[K_n:K\right]\dd_{\L,K_n}(P)\right)^{2g+2\kappa(g)}.\]
\end{theo}

\medskip

\rem Un peu plus g\'en\'eralement le m\^eme \'enonc\'e est valable pour toute extension $L/K$ ab\'elienne telle que son discriminant $\text{disc}(L/K)$ v\'erifie 
\begin{equation}\label{ra}
\frac{1}{[L:K]}\log \text{disc}(L/K)\leq c_1(A/K,\mathcal{L})\left(\log[L:K]\right)^2.
\end{equation}
\noindent Notamment (cf. la preuve du lemme \ref{rami}) ceci est vrai pour toute extension $K\subset L\subset K_n$ telle que $\log\phi(n)\ll \log [L:K]$. Par exemple, on peut voir que le th\'eor\`eme pr\'ec\'edent reste valable pour tout $n$ avec $L=K(T_n)$ o\`u $T_n$ est un point de torsion de $A$ d'ordre $n$.\footnote{Notons que l'in\'egalit\'e (\ref{ra}) n'est pas toujours v\'erifi\'ee comme le montre l'exemple $\Q\subset \Q(\sqrt{l})\subset \Q(A[l])$ avec $l$ premier congru \`a $1$ modulo $4$.}

\medskip

\noindent Comme dit pr\'ec\'edemment il s'agit ici de l'analogue d'un r\'esultat de Amoroso et David obtenu dans \cite{ad2003} pour le groupe multiplicatif $\mathbb{G}_m^{n}$. 

\medskip

\rem En appliquant la remarque 1.2, on d\'eduit du th\'eor\`eme \ref{semi} le m\^eme \'enonc\'e avec $\dd_{\L,K_n}(P)$ remplac\'e par le degr\'e $[K_n(P):K_n]^{\frac{1}{g}}$.

\medskip

\noindent On a ici quantifi\'e l'hypoth\`ese ``$P$ est d'ordre infini'' par une borne sur le degr\'e de la plus petite sous-vari\'et\'e de torsion contenant le point $P$. Notons que ce type de r\'esultats, rarement mis en \'evidence, contient la plupart du temps des informations arithm\'etiques suppl\'ementaires. On donne d'ailleurs ici, suivant la preuve du r\'esultat analogue pour $\G_m^n$ due \`a Amoroso-David \cite{ad2003}, une cons\'equence de cette quantification. Il s'agit d'une minoration du minimum absolu des sous-vari\'et\'es non de torsion de $A^n$ o\`u $A/K$ est une vari\'et\'e ab\'elienne simple de type C.M. Pour cela on introduit deux notations : on pose $\L_n$ le fibr\'e $\L^{\boxtimes n}$ sur $A^n$ d\'eduit d'un fibr\'e ample et sym\'etrique $\L$ sur $A$. Par ailleurs, on note 
\[V^{\star}=V\backslash\bigcup B,\]
\noindent o\`u l'union porte sur les sous-vari\'et\'es $B$ de torsion incluses dans $V$. On donne dans l'\'enonc\'e suivant une minoration du minimum absolu d'une vari\'et\'e non de torsion. Notons que l'on ne la suppose pas irr\'eductible.

\medskip

\defi Supposons que la vari\'et\'e ab\'elienne $A$ soit plong\'ee dans un espace projectif $\mathbb{P}_n$; notons $\mathcal{R}$ l'anneau des
coordonn\'ees de $\mathbb{P}_n$; on dira qu'une sous-vari\'et\'e $V$ de $A$ peut \^etre \textit{d\'efinie incompl\`etement
par des \'equations de degr\'e $\leq L$} si $V$ est une composante isol\'ee de $A\cap Z(\mathcal{I})$  o\`u $\mathcal{I}$ est
un id\'eal homog\`ene de $\mathcal{R}$ engendr\'e par des polyn\^omes de degr\'e au plus $L$.

\begin{theo}\label{absolu} Soit $A/K$ une vari\'et\'e ab\'elienne simple de type C.M., plong\'ee dans un espace projectif par un fibr\'e en droites ample $\L$, et $n\geq 1$ un entier. Il existe une constante strictement positive $c(A/K,\L,n)$ telle que : pour toute sous-vari\'et\'e (non n\'ecessai\-re\-ment irr\'eductible) $V$ de $A^n$, d\'efinie sur une extension $K_r:=K(A[r])/K$ pour un certain entier $r$, incompl\`etement d\'efinie dans $A^n$ sur $K_r$ par des \'equations de degr\'e au plus $\ddd$, alors on a pour tout point $\overline{K}$-rationnel $P$ de $V^{\star}$~:
\[ \h_{\L_n}(P)\geq \frac{c(A/K,\L,n)}{\ddd}\left(\frac{\log\log 3\left[K_r:K\right]\ddd}{\log 3\left[K_r:K\right]\ddd}\right)^{\kappa(g)},\]
\noindent avec  $\kappa(g)=(g+1)!(2g+5)(g+2)(2g.g!)^g$.
\end{theo}

\medskip

\noindent Notons (\textit{cf}. par exemple \cite{bombierizannier} p.790) que $V$ est toujours incompl\`etement d\'efinie par des \'equations de degr\'e au plus $c_1(A/K,\mathcal{L},n)\deg_{\L} V$. En notant $\mu^{\star}(V)$ le minimum absolu de $V$ (\textit{i.e.} le minimum des hauteurs des points de $V^{\star}$), on obtient en corollaire :

\medskip

\begin{cor}\label{corabsolu} Soit $A/K$ une vari\'et\'e ab\'elienne simple de type C.M. et $n\geq 1$ un entier. Il existe une constante strictement positive $c(A/K,n)$ telle que : pour toute sous-vari\'et\'e (non n\'ecessai\-re\-ment irr\'eductible) $V/K$ de $A^n$, non de torsion et d\'efinie sur $K$, on a:
\[ \mu^{\star}(V)\geq \frac{c(A/K,n)}{\deg_{\L}V}\left(\frac{\log\log 3\deg_{\L}V}{\log 3\deg_{\L}V}\right)^{\kappa(g)},\]
\noindent avec  $\kappa(g)=(g+1)!(2g+5)(g+2)(2g.g!)^g$.
\end{cor}

\medskip

\noindent Ceci rend effectif les r\'esultats de Bombieri et Zannier \cite{bombierizannier}, et raffine ceux de David et Philippon \cite{dp} concernant ce m\^eme probl\`eme. Plus exactement en suivant la preuve de \cite{bombierizannier}, on constate qu'ils obtiennent une minoration effective mais qui est multiexponentielle en le degr\'e, au lieu d'\^etre comme ici lin\'eaire. Concernant le r\'esultat de \cite{dp}, bien qu'il ne soit pas explicitement r\'edig\'e, on peut voir qu'ils obtiennent comme cons\'equence de leur th\'eor\`eme principal, une minoration polynomiale en le degr\'e.

\medskip

\subsubsection{Intersection de courbes et de sous-groupes}

\noindent Soient $A/\overline{\Q}$ une vari\'et\'e ab\'elienne (d\'efinie) sur $\overline{\Q}$, $X/\overline{\Q}$ une courbe g\'eom\'etri\-que\-ment irr\'educ\-tible de $A$ et $r$ un entier positif ou nul. Comme dit pr\'ec\'edemment on consid\`ere l'ensemble
\[A^{[r]}:=\bigcup_{\text{codim\,}G\geq r}G(\overline{\Q})\]
\noindent o\`u l'union porte sur les sous-groupes alg\'ebriques non-n\'ecessairement connexes de $A$ de codimension au moins $r$. De mani\`ere ind\'ependante, Zilber (\cite{zilber} Conjecture 2) pour les vari\'et\'es semi-ab\'eliennes et Pink (\cite{pink} Conjecture 1.3) pour les vari\'et\'es de Shimura mixtes, ont formul\'es une conjecture qui, dans le cas des courbes incluses dans une vari\'et\'e ab\'elienne sur $\overline{\Q}$ se sp\'ecialise en la suivante : 

\medskip

\begin{conj}\label{conji}\textnormal{\textbf{(Zilber-Pink, cas particulier)}} Soient $A/\overline{\Q}$ une vari\'et\'e ab\'elienne et $X/\overline{\Q}$ une courbe dans $A$ g\'eom\'etriquement irr\'eductible. Si $X$ n'est pas contenue dans un sous-groupe alg\'ebrique strict de $A$, alors l'ensemble $X(\overline{\Q})\cap A^{[2]}$ est fini.
\end{conj}

\medskip

\noindent R\'emond \cite{remond} a montr\'e que la conjecture pr\'ec\'edente est vraie si une tr\`es bonne minoration (conjecturale) des points d'ordre infini de $A$ est vraie~: 

\medskip

\begin{theor}\label{condi}\textnormal{\textbf{(R\'emond \cite{remond})}} Si la conjecture \ref{conj2} est vraie pour toute vari\'et\'e ab\'elienne alors la conjecture \ref{conji} est vraie pour toute vari\'et\'e ab\'elienne.
\end{theor}

\medskip

\noindent Par ailleurs dans le cas des vari\'et\'es ab\'eliennes de type CM, R\'emond \cite{remond} obtient un r\'esultat inconditionnel, mais sensiblement plus faible que la conjecture \ref{conji}~: soit $A$ une vari\'et\'e ab\'elienne de type CM, g\'eom\'etriquement isog\`ene au produit $\prod_{i=1}^mA_i^{n_i}$ o\`u les $A_i$ sont des vari\'et\'es ab\'eliennes simples de dimension respective $g_i$, deux \`a deux non-isog\`enes.

\medskip

\begin{theor}\textnormal{\textbf{(R\'emond \cite{remond})}}\label{gr} Soient $A/\overline{\Q}$ une vari\'et\'e ab\'elienne de type CM et $X$ une courbe dans $A$ non contenue dans un sous-groupe alg\'ebrique strict, alors $X(\overline{\Q})\cap A^{[2+\sum_{i=1}^m g_i]}$ est fini.
\end{theor}

\noindent En utilisant notre r\'esultat de minoration de hauteur (\textit{cf.} th\'eor\`eme \ref{semi} ci-dessus) ainsi qu'un r\'esultat faisant l'objet d'un article s\'epar\'e \cite{matorsion} concernant la torsion dans les vari\'et\'es ab\'eliennes de type C.M., nous am\'eliorons ceci et obtenons un r\'esultat optimal dans le cas d'une puissance d'une vari\'et\'e ab\'elienne simple de type C.M.~:

\medskip

\begin{theo}\label{ray} La conjecture \ref{conji} est vraie pour toute vari\'et\'e ab\'elienne $A$ de type C.M., isog\`ene \`a une puissance d'une vari\'et\'e ab\'elienne simple. 
\end{theo}

\medskip

\rem Dans les \'enonc\'es pr\'ec\'edents (th\'eor\`emes \ref{gr} et \ref{ray}) nous avons utilis\'e l'hypoth\`ese \textit{$X$ contenue dans aucun sous-sch\'ema en groupes de $A$, distinct de $A$}. Dans son article \cite{remond}, R\'emond utilise en fait une hypoth\`ese plus forte : il suppose que $X$ est transverse (\textit{i.e.} contenue dans aucune translat\'ee de sous-vari\'et\'e ab\'elienne). En fait le corollaire 1.1 de \cite{remond2} et sa preuve montrent que l'on peut dans nos \'enonc\'es dans le cas des vari\'et\'es ab\'eliennes,  remplacer l'hypoth\`ese \textit{$X$ transverse} par l'hypoth\`ese plus faible \textit{$X$ contenue dans aucun sous-sch\'ema en groupes de $A$, distinct de $A$}. C'est ce que nous avons fait ici.  

\medskip

\noindent On peut en fait \'enoncer un r\'esultat optimal dans un cadre un peu plus vaste qu'une puissance d'une vari\'et\'e ab\'elienne simple de type C.M. Consid\'erons pour cela une notion introduite formellement dans \cite{matorsion}~:

\medskip

\defi Soit $A/K$ une vari\'et\'e ab\'elien\-ne quelconque sur un corps de nombres. On d\'efinit un invariant $\gamma(A)$ par
\[\gamma(A)=\inf\left\lbrace x>0\ | \ \exists C>0,\ \ \forall F/K \text{ finie, }\ \left|A(F)_{\tors}\right|\leq C[F:K]^x\right\rbrace.\]

\noindent Le seul r\'esultat connu en toute g\'en\'eralit\'e pour cet invariant est d\^u \`a Masser \cite{lettre}~: il obtient 
\[\gamma(A)\leq g.\]
\noindent Par ailleurs, dans \cite{matorsion} l'auteur a obtenu une reformulation combinatoire de cet invariant $\gamma(A)$ dans le cas d'une vari\'et\'e ab\'elienne $A/K$ g\'eom\'etriquement simple de type C.M. En particulier, dans le cas d'une telle vari\'et\'e ab\'elienne, le corollaire 1.13 de \cite{matorsion} donne la majoration
\[\gamma(A)\leq\frac{2g}{2+\log_2g}\]
\noindent o\`u $g=\dim A$ et $\log_2$ est le $\log$ en base $2$. Cette majoration raffine celle de Masser~:  si $g\geq 2$ alors $\gamma(A)<g$. De plus dans le cas d'une vari\'et\'e ab\'elienne CM de type non-d\'eg\'en\'er\'e (\textit{i.e.} ayant un groupe de Mumford-Tate aussi grand que possible, donc de dimension $g+1$) on montre m\^eme que $\gamma(A)$ est beaucoup plus petit : $\gamma(A)=\frac{2g}{g+1}$.

\medskip

\noindent Soit $A/K$ une vari\'et\'e ab\'elienne. Elle est isog\`ene \`a un produit :  $\prod_{i=1}^mA_i^{n_i}$, o\`u les $A_i$ sont simples et deux \`a deux non-isog\`enes, de dimension $g_i$. On montre au paragraphe \ref{pararay} que si $X$ est une courbe transverse de $A$ et si $A$ est de type C.M., alors
\[\min_{1\leq i\leq m}g_i>\sum_{i=1}^m\gamma(A_i)\Longrightarrow \left|X(\overline{K})\cap A^{[2]}\right|<+\infty.\]

\noindent \textbf{Plan de l'article : }les paragraphes 2 et 3 sont consacr\'es \`a la preuve du th\'eor\`eme \ref{semi} concernant le probl\`eme de Lehmer, le paragraphe 4 est consacr\'e au th\'eor\`eme \ref{ray} concernant le probl\`eme d'intersection de courbe et de sous-groupes et le paragraphe 5 au th\'eor\`eme \ref{absolu} concernant le minimum absolu. L'appendice explique les liens entre les parties (1) et (2) du probl\`eme de Lehmer ab\'elien ainsi que les liens entre le probl\`eme de Lehmer ab\'elien et sa variante 
multihomog\`ene (conjecture \ref{conj2app2}) formul\'ee par David et Hindry.

\medskip

\noindent Plus pr\'ecis\'ement, au paragraphe 2 nous donnons une majoration du discriminant absolu de l'extension $K(A[n])/\mathbb{Q}$ dont nous aurons besoin dans l'application du lemme de Siegel. Nous faisons \'egalement dans ce paragraphe les rappels n\'ecessaires sur les isog\'enies de Frobenius. Classiquement depuis le r\'esultat de Dobrowolski concernant le probl\`eme de Lehmer, c'est essentiellement en utilisant des transform\'es par ces isog\'enies d'un point $P$, dont on suppose par l'absurde qu'il ne v\'erifie pas la conclusion du th\'eor\`eme \ref{semi}, que nous allons extrapoler. Il y a ceci dit ici une diff\'erence, dans la mesure o\`u l'on travaille, non plus sur $K$ mais sur $K(A[n])$ : on consid\'erera plut\^ot des tordus de ces points par un certain automorphisme de $\Gal(\overline{K}/K)$, un Frobenius de l'extension ab\'elienne $K(A[n])/K$. Au paragraphe 3, et avec la diff\'erence pr\'ec\'edemment indiqu\'ee, nous donnons la preuve du th\'eor\`eme \ref{semi} en suivant ce qui est fait dans \cite{davidhindry} et avec les id\'ees nouvelles introduites dans \cite{ad2003}. Notons tout de m\^eme une autre diff\'erence par rapport \`a \cite{davidhindry}~: on travaille avec les indices d'obstruction alors que dans \cite{davidhindry} la preuve est faite avec le degr\'e $D$. Il y a donc certaines modifications suppl\'ementaires \`a faire. Ceci \'etant la preuve est une classique preuve d'approximation diophantienne, avec un lemme de Siegel (le th\'eor\`eme de Bombieri-Vaaler qui nous permet d'avoir un contr\^ole explicite de la d\'ependance en le corps sur lequel on travaille), un lemme de z\'eros et, comme dans \cite{davidhindry} une descente finale o\`u l'on r\'eit\`ere $g$ fois l'ensemble de la preuve pour pouvoir conclure. Dans ce paragraphe 3 nous avons choisi d'insister plus particuli\`erement sur les points qui diff\`erent par rapport \`a \cite{davidhindry} en renvoyant \`a cette r\'ef\'erence lorsque cela s'av\'erait n\'ecessaire.

\noindent Au paragraphe 4, nous donnons une preuve, suivant la strat\'egie de \cite{remond}, du th\'eor\`eme \ref{ray} en utilisant notre r\'esultat sur le probl\`eme de Lehmer.

\noindent Au paragraphe 5 nous donnons une preuve par r\'ecurrence sur $n$ du th\'eor\`eme \ref{absolu} bas\'e sur notre th\'eor\`eme \ref{semi} et suivant la preuve de \cite{ad2003} dans le cas de $\G_m^n$. On utilise ici \`a plein l'alternative d\'emontr\'ee dans le th\'eor\`eme \ref{semi}.

\noindent Enfin dans l'appendice nous montrons, que la partie (2) du probl\`eme de Lehmer ab\'elien \ref{conj1app2} entra\^ine sa partie (1) ; autrement dit, si l'on sait minorer la hauteur des points qui ne sont contenus dans aucune sous-vari\'et\'e de torsion, alors on sait en fait aussi bien minorer la hauteur de tous les points non de torsion. Ceci nous permet de raffiner, dans le cas des vari\'et\'es ab\'eliennes de type C.M. un r\'esultat de Masser \cite{lettre} sur la minoration des points d'ordre infini. Dans le second et dernier sous-paragraphe de l'appendice, nous montrons que la conjecture multihomog\`ene de David et Hindry \ref{conj2app2} est en fait une cons\'equence du probl\`eme de Lehmer ab\'elien.

\section{Pr\'eliminaires arithm\'etico-g\'eom\'etriques}

\noindent Dans la suite de l'article nous introduisons un certain nombre de constantes $c_i$. Nous indiquerons en g\'en\'eral entre parenth\`eses de quoi depend cette constante (par exemple $c_i(A/K)$ depend de $A/K$, $c_i(\L)$ d\'epend de $\L$...). Si toutefois aucune d\'ependance n'est explicitement indiqu\'ee, c'est que la constante d\'epend uniquement des donn\'ees $A/K$ et $\L$.

\subsection{Ramification}

\noindent On travaille ici et comme partout avec une vari\'et\'e ab\'elienne $A/K$ o\`u $K$ est un corps de nombres.

\medskip

\noindent \textbf{Notations.} On note $\phi$ la \textit{fonction d'Euler} et, si $L/F$ est une extension finie de corps, on note $\Dif_{L/F}$ la \textit{diff\'erente de $L/F$}. Par ailleurs, si $n$ est un entier, on note $K_n=K(A[n])$ l'extension de degr\'e $\deg(K_n/K)$ de $K$ engendr\'ee par les points de torsion d'ordre divisant $n$ de $A(\overline{K})$. Dans ce qui suit, on s'int\'eresse \`a une majoration du \textit{discriminant absolu}, $\text{disc}(K_n)$, de $K_n/\Q$. On donne de plus une majoration du nombre et de la taille des premiers qui se ramifient dans $K_n$.

\medskip

\begin{lemme}\label{nombrepremier}Si $n$ un entier strictement positif, on note $\omega(n)$ le nombre de premiers deux \`a deux distincts divisant $n$. On a les in\'egalit\'es :
\begin{enumerate}
\item Pour tout entier $n\geq 2$, on a $\omega(n)\leq\frac{2\log\phi(n)}{\log 2}.$
\item Pour tout entier $n$ assez grand, on a $\omega(n)\leq\frac{4\log n}{\log\log n}.$
\end{enumerate}
\end{lemme}
\demo Le point 1. est un exercice facile et le point 2. est un r\'esultat classique que l'on peut par exemple trouver dans \cite{tenenbaum}.\hfill$\Box$

\medskip

\begin{lemme}\label{nombrepremier2}Soit $n$ un entier strictement positif. On \'ecrit sa d\'ecomposition en facteurs premiers $n=\prod p_i^{\alpha_i}$. On a 
\[\prod_{i=1}^{\omega(n)}p_i\leq \phi(n)^3.\]
\end{lemme}
\demo C'est un simple calcul :
\begin{align*}
\prod_{i=1}^{\omega(n)} p_i	& \leq \prod_{\alpha_i\geq 2}\phi(p_i^{\alpha_i})\prod_{\alpha_i=1}(\phi(p_i)+1)\leq 2^{\text{Card}\left\{i | \alpha_i=1\right\}}\prod_{\alpha_i\geq 2}\phi(p_i^{\alpha_i})\prod_{\alpha_i=1}\phi(p_i)\\
		& \leq  2^{\text{Card}\left\{p | p\mid n\right\}}\phi(n)\leq 2^{\frac{2\log\phi(n)}{\log 2}}\phi(n)=\phi(n)^3,\\
\end{align*}
\noindent o\`u l'on a utilis\'e le lemme \ref{nombrepremier} pr\'ec\'edent dans la derni\`ere in\'egalit\'e.\hfill$\Box$

\medskip

\begin{lemme}\label{rami} On suppose ici que $A/K$ a partout bonne r\'eduction. 
\begin{enumerate}
\item Il existe une constante strictement positive $c_1([K:\Q])$ telle que pour tout $n\geq 1$, on a 
\[ \text{disc}(K_n)\leq \left(\deg(K_n/K)\right)^{c_1\deg(K_n/K)\log \deg(K_n/K)}.\]
\item Il existe une constante $C_1(K)$ telle que le cardinal des premiers ramifi\'es dans $K_n$ est major\'e par $C_1(K)+\log\phi(n)$.
\end{enumerate}
\end{lemme}
\demo Pour le point (1) : soient $\p$ un id\'eal premier de $K$ au-dessus d'un nombre premier $p$ et $\pp$ un id\'eal premier de $K_n$ au-dessus de $\p$. Soient $v$ la place associ\'ee \`a $p$, $w_K$ la place associ\'ee \`a $\p$ et $w_{K_n}$ celle associ\'ee \`a $\pp$. Soient $\pi$, $\pi_K$, $\pi_{K_n}$ les uniformisantes associ\'ees \`a ces id\'eaux. Si $e(\p)$ (respectivement $e(\pp)$) est l'indice de ramification de $\p$ sur $p$ (respectivement de $\pp$ sur $\p$), on a
\[w_{K_n}(\pi)=e(\p)w_{K_n}(\pi_K)=e(\p)e(\pp)w_{K_n}(\pi_{K_n})=e(\p)e(\pp), \text{ et, } w_{K}(\pi)=e(\p).\]
\noindent On sait par le Theorem 1. de \cite{serretate} (voir \'egalement la proposition 18 de \cite{shita}) que, si la vari\'et\'e ab\'elienne $A$ a bonne r\'eduction en $\p$, alors l'extension $K_n/K$ ne peut \^etre ramifi\'ee en $\p$ que si $p$ divise $n$ (c'est le sens facile du crit\`ere de N\'eron-Ogg-Schafarevitch). 
 On a donc
\[\Dif_{K_n/K}=\prod_{\pp|p\mid n}\pp^{\ord_{\pp}(\Dif_{K_n/K})}.\]
\noindent Or la proposition 13 de \cite{corpslocaux} et la remarque suivant cette proposition nous donne la borne :
\[\ord_{\pp}(\Dif_{K_n/K})\leq e(\pp)-1+w_{K_n}(e(\pp))=e(\pp)-1+e(\p)e(\pp)v(e(\pp)).\]
\noindent Par ailleurs, l'indice $e(\pp)$ \'etant inf\'erieur \`a $\deg(K_n/K)$, on a $v(e(\pp))\leq \frac{\log \deg(K_n/K)}{\log 2}.$ Ainsi, on obtient la minoration 
\[\ord_{\pp}(\Dif_{K_n/K})\leq \frac{2}{\log 2}e(\pp)\log \deg(K_n/K).\]
Par transitivit\'e et par d\'efinition du discriminant $\text{disc}(K_n)$ on a 
\begin{align*}
\text{disc}(K_n)	& =N_{\Q}^K\left(\text{Disc}(K_n/K)\right)\text{Disc}(K/\Q)^{\deg(K_n/K)} = c_2^{\deg(K_n/K)}N_{\Q}^{K_n}\left(\Dif(K_n/K)\right),\\
	& \leq \prod_{p\mid n} p^{c_3\sum_{\pp|p}f(\pp)e(\pp)\log \deg(K_n/K)} \leq \prod_{p\mid n} p^{c_3\deg(K_n/K)\log \deg(K_n/K)}.\\
\end{align*}
\noindent Le formalisme de l'accouplement de Weil nous indique que $\mu_n\subset K_n$, donc pour les degr\'es que $\deg(K_n/\mathbb{Q})\geq \phi(n)$. Ainsi, en utilisant le lemme \ref{nombrepremier2} pr\'ec\'edent, on peut conclure :
\[\text{disc}(K_n)\leq \phi(n)^{3c_3\deg(K_n/K)\log \deg(K_n/K)}\leq \deg(K_n/\mathbb{Q})^{3c_3\deg(K_n/K)\log \deg(K_n/K)}.\]
\noindent Concernant le point (2), on constate que $p$ est ramifi\'e dans $K_n$ si et seulement si il est ramifi\'e dans $K$ ou s'il est non ramifi\'e dans $K$ mais tel qu'il existe un premier $\p$ de $K$ au dessus de $p$ se ramifiant dans $K_n$. Le sens facile du crit\`ere de N\'eron-Ogg-Schafarevitch et le point (1) du lemme \ref{nombrepremier} permettent de conclure.\hfill$\Box$

\medskip

\begin{cor}\label{coro}Il existe une constante strictement positive $c_2$ telle que pour tout $n\geq 1$, on a 
\[ \log \left(\text{disc}(K_n)^{\frac{1}{\deg(K_n/K)}}\right)\leq c_2\left(\log \deg(K_n/K)\right)^2.\]
\end{cor} 
\demo C'est imm\'ediat.\hfill $\Box$

\medskip

\rem Cette estimation n'est tr\`es certainement pas optimale, mais elle est simple \`a obtenir et nous suffira.

\medskip

\subsection{Frobenius et isog\'enies admissibles}

\subsubsection{Morphismes de Frobenius}

\noindent \textbf{Notations.} Si $K$ est un corps de nombres, on note $\mathcal{O}_K$ son anneau d'entiers, $v$ une place finie de $K$, et $k_v$ le corps r\'esiduel associ\'e \`a $v$.

\medskip

\noindent Si $A/K$ est une vari\'et\'e ab\'elienne, on note $\mathcal{A}/\mathcal{O}_K$ son mod\`ele de N\'eron, et $A_v/k_v$ la fibre sp\'eciale correspondant \`a la place finie $v$. Rappelons la propri\'et\'e universelle du mod\`ele de N\'eron : si $\mathcal{X}/\mathcal{O}_K$ est lisse, de fibre g\'en\'erique $X/K$, tout $K$-morphisme $X\rightarrow A$ se rel\`eve de mani\`ere unique en un $\mathcal{O}_K$-morphisme $\mathcal{X}\rightarrow \mathcal{A}$. La propri\'et\'e universelle du produit fibr\'e $A_v=\mathcal{A}\times_{\mathcal{O}_K}k_v$ permet d'associer naturellement \`a tout $\mathcal{O}_K$-endomorphisme de $\mathcal{A}$ un $k_v$-endomorphisme de $A_v$. En utilisant la propri\'et\'e universelle du mod\`ele de N\'eron, on en d\'eduit une fl\`eche naturelle 
\[\Psi : \textnormal{End}_K(A) \rightarrow \textnormal{End}_{k_v}(A_v).\]
Cette fl\`eche n'est en g\'en\'eral pas surjective, mais on peut par contre montrer qu'elle est injective aux places de bonne r\'eduction.

\medskip

\noindent Sur la vari\'et\'e $A_v/k_v$, on dispose d'un endomorphisme particulier : le morphisme de Frobenius $\textnormal{Frob}_v$, correspondant en coordonn\'ees projectives \`a l'\'el\'e\-va\-tion \`a la puissance $q=\textnormal{N}(v)$, o\`u $\textnormal{N}(v)$ est la norme $K/\mathbb{Q}$ de $v$. Dans le cas C.M., un th\'eor\`eme de Shimura-Taniyama permet d'affirmer que le morphisme $\textnormal{Frob}_v$ se rel\`eve en presque toute place :

\begin{propo}\label{shimura}\textnormal{\textbf{(Shimura-Taniyama)}} Soit $A/K$ une vari\'et\'e ab\'elienne de ty\-pe C.M. Notons $\prod_{i=1}^r K_i$ le produit de corps de nombres qui est inclus dans $\textnormal{End}_{\overline{K}}(A)\otimes\mathbb{Q}$ et tel que $\sum_{i=1}^r[K_i :\mathbb{Q}]=2\,\textnormal{dim} A$. On suppose que le corps de nombres $K$ contient tous les $K_i$, et que $\prod_{i=1}^r\mathcal{O}_{K_i}$ est inclus dans $\textnormal{End}_K(A)$. Alors, pour toutes les places sauf \'evntuellement un nombre fini d'entre elles, l'endomorphisme $\textnormal{Frob}_v$ se rel\`eve en un $K$-endomorphisme $\alpha_v$ de $A$. On appelera morphisme de Frobenius sur $A$ un tel endomorphisme.
\end{propo}
\noindent \textbf{D\'emonstration} C'est le Theorem 1 paragraphe III.13 de \cite{shita}.\hfill$\Box$

\medskip

\noindent Ce sont ces morphismes de Frobenius sur $A/K$ qui vont nous permettre d'\'ecrire l'\'etape d'extrapolation.

\medskip

\noindent Enfin pour tout nombre premier $p$, on note $\Phi_p$ l'automorphisme de Frobenius de l'extension ab\'elienne $K_n/K$ et on choisit une extension de $\Phi_p$ \`a $\Gal(\overline{K}/K)$, que l'on note encore $\Phi_p$. 

\subsubsection{Isog\'enies admissibles}
\noindent On rappelle la notion d'isog\'enie admissible telle qu'introduite dans \cite{davidhindry}.

\medskip

\defi Soient $A$ une vari\'et\'e ab\'elienne et $\mathcal{L}$ un fibr\'e en droites ample sur $A$. Une isog\'enie $\alpha$ de $A$ est dite \textit{admissible} par rapport \`a $\mathcal{L}$ si 
\begin{enumerate}
\item $\alpha$ est dans le centre de $\textnormal{End}(A)$.

\medskip

\item il existe un entier $\textnormal{q}(\alpha)$ appel\'e \textit{poids} de $\alpha$ tel que $\alpha^{\star}\mathcal{L}\simeq \mathcal{L}^{\otimes \, \textnormal{q}(\alpha)}$.
\end{enumerate}

\medskip

\rem En fait la condition (1) ne sert qu'a simplifier l'\'enonc\'e du lemme \ref{distinct}. C'est la condition (2) qui importe vraiment. Les seules isog\'enies qui nous int\'eresseront sont les relev\'ees $\alpha_v$ des morphismes de Frobenius qui sont admissibles (\textit{cf.} la Proposition \ref{frob}).

\begin{lem}\label{degre}\textnormal{\textbf{(David-Hindry)}} Soient $A$ une vari\'et\'e ab\'elienne de dimension $g$ munie d'un fibr\'e en droites tr\`es ample $\mathcal{L}$, et $\alpha$ une isog\'enie admissible relativement \`a $\mathcal{L}$, de poids $q=\textnormal{q}(\alpha)$. Dans le plongement projectif de $A$, associ\'e \`a $\mathcal{L}$, $A\hookrightarrow \mathbb{P}_n$, on a : 
\begin{enumerate}
\item $\textnormal{card}\left(\textnormal{Ker}(\alpha)\right)=q^g$,
\item pour toute sous-vari\'et\'e $V$ de $A$ de stabilisateur $G_V$, on a 
\[\deg_{\mathcal{L}}\left(\alpha(V)\right)=\frac{q^{\dim(V)}}{\left| G_V\cap \textnormal{Ker}(\alpha)\right|}\deg_{\mathcal{L}}(V)\]
\item pour toute sous-vari\'et\'e $V$ de $A$, d\'efinie incompl\`etement dans $A$ par des \'equations de degr\'e inf\'erieur \`a $L$, de stabilisateur $G_V$ on a 
\[\deg_{\L} G_V=\mid G_V:G_V^0\mid\deg_{\L} G_V^0\leq \deg_{\L} (V)(2L)^{\dim V-\dim G_V}\]
\noindent et $G_V$ est d\'efini incompl\`etement dans $A$ par des \'equations de degr\'e inf\'erieur \`a $2L$.
\end{enumerate}
\end{lem}
\noindent \textbf{D\'emonstration} Le point (1) est facile : par d\'efinition, $\alpha^{\star}\mathcal{L}\simeq \mathcal{L}^{\otimes q}$. On a donc, 
\[q^g\deg_{\mathcal{L}}(A)=\deg_{\mathcal{L}^{\otimes q}}(A)=\deg_{\alpha^{\star}\mathcal{L}}(A)=\left|\textnormal{Ker}(\alpha)\right|\deg_{\mathcal{L}}(A).\]
\noindent L'amplitude de $\mathcal{L}$ nous assure que le dernier degr\'e est strictement positif. On simplifie pour conclure. Pour le point (2), il s'agit du point (ii) du lemme 6. de \cite{hindry} et le point 3. correspond au point (ii) du lemme 2.1. de \cite{davidhindry}.\hfill$\Box$

\medskip

\begin{lemme}\label{lemme25}Soient $G$ un sous-groupe alg\'ebrique de la vari\'et\'e ab\'elienne $A/K$, $\mathcal{L}$ un fibr\'e en droites tr\`es ample sur $A$, et $\alpha$ une isog\'enie admissible relativement \`a $\mathcal{L}$ de poids $\textnormal{q}(\alpha)$ de $A$. On a 
\[ \textnormal{q}(\alpha)^{\dim G}\leq\textnormal{card}\left(\textnormal{Ker}(\alpha)\cap G\right)\leq \left[G:G^0\right]\textnormal{q}(\alpha)^{\dim G}.\]
\end{lemme}
\noindent \textbf{D\'emonstration} On note que 
\[ \left[G:G^0\right]\textnormal{card}\left(\textnormal{Ker}(\alpha)\cap G^0\right)\geq \textnormal{card}\left(\textnormal{Ker}(\alpha)\cap G\right)\geq \textnormal{card}\left(\textnormal{Ker}(\alpha)\cap G^0\right).\]
\noindent La restriction de $\alpha$ \`a la sous-vari\'et\'e ab\'elienne $G^0$ est encore une isog\'enie admissible de poids $\textnormal{q}(\alpha)$ pour $(G^0,\mathcal{L}_{\mid G^0})$ (\textit{cf.} Lemme 2.4. point (ii) de \cite{davidhindry}). Par le point (1) du lemme \ref{degre} pr\'ec\'edent, on en d\'eduit que le cardinal du noyau de cette isog\'enie $\alpha_{\mid G^0}$ est $\textnormal{q}(\alpha)^{\dim G^0}$. \hfill $\Box$

\medskip

\noindent Soient $F/K$ une extension finie de corps et $V$ une sous-$F$-vari\'et\'e stricte de $A_F$ (produit fibr\'e de $A$ par Spec $F$ au dessus de Spec $K$), $F$-irr\'educ\-tible. Le lemme suivant (dont l'origine remonte \`a Dobrowolski \cite{dob}) montre que les images par une isog\'enie admissible de ses composantes g\'eom\'etriquement irr\'eductibles sont essentiellement distinctes. Nous en aurons besoin au paragraphe \ref{fin}. On commence pour cela par donner une d\'efinition~:

\medskip

\defi Soient $A$ une vari\'et\'e ab\'elienne et $\mathcal{L}$ un fibr\'e en droites ample sur $A$. Deux isog\'enies admissibles de $A$ par rapport \`a $\mathcal{L}$ sont dites \textit{premi\`eres entre elles} si leurs poids sont premiers entre eux.

\begin{lemme}\label{distinct} Soient $A$ une vari\'et\'e ab\'elienne sur $K$ de dimension $g\geq 1$, $\mathcal{L}$ un fibr\'e en droites tr\`es ample sur $A$, $F/K$ une extension finie de corps de nombres et $V$ une sous-$F$-vari\'et\'e stricte de $A_F$, irr\'eductible sur $F$.  Si $V_{\overline{K}}$ n'est pas une r\'eunion de sous-vari\'et\'es de torsion de $A_{\overline{K}}$, on a :
\begin{enumerate}
\item Pour tout couple $(\alpha,\beta)$ d'isog\'enies admissibles pour $\mathcal{L}$, de poids distincts, pour tout $\sigma\in \textnormal{Gal}(\overline{K}/F)$, et pour toute composante g\'eom\'etri\-quement irr\'e\-ductible $W$ de $V_{\overline{K}}$, les sous-vari\'et\'es $\alpha(W)$ et $\beta\left(\sigma(W)\right)$ sont distinctes.

\medskip

\item Soit $\mathcal{P}$ un ensemble d'isog\'enies admissibles pour $\mathcal{L}$, deux \`a deux pre\-mi\`eres entre elles. Notons $V_1,\ldots, V_M$ les composantes g\'eom\'etri\-que\-ment irr\'e\-duc\-tibles de $V_{\overline{K}}$, et notons $\mathcal{Q}$ le sous-ensemble de $\mathcal{P}$ d\'efini par
\[ \mathcal{Q}=\left\{ \alpha\in \mathcal{P}\ | \ \exists i,j,\ 1\leq i<j\leq M,\ \ \alpha(V_i)=\alpha(V_j)\right\}.\]
\noindent Le cardinal de $\mathcal{Q}$ est major\'e par $\frac{\log M}{\log 2}$.
\end{enumerate}
\end{lemme}
\noindent \textbf{D\'emonstration} Dans ce contexte il s'agit de la proposition 2.7. de \cite{davidhindry} appliqu\'ee sur le corps de nombres $F$.\hfill$\Box$

\medskip

\noindent On conclut ce paragraphe en ``rappelant'' que les morphismes de Frobenius sur $A/K$ sont des isog\'enies admissibles :

\medskip

\defi Soient $A$ une vari\'et\'e ab\'elienne et $\mathcal{L}$ un fibr\'e en droites ample sur $A$. Suivant Mumford, on dit que $\mathcal{L}$ est \textit{totalement sym\'etrique} si $\mathcal{L}$ est le carr\'e d'un fibr\'e en droites sym\'etrique.

\medskip

\noindent Le th\'eor\`eme de Lefschetz (\textit{cf.} par exemple le Theorem A.5.3.6 de \cite{hindrysil}) nous indique que si $\mathcal{L}$ est un fibr\'e en droites ample, alors $\mathcal{L}^{\otimes 3}$ est tr\`es ample.

\begin{prop}\label{frob}Soient $A/K$ une vari\'et\'e ab\'elienne de type C.M. v\'erifiant les hypoth\`eses de la proposition \ref{shimura}, et $\mathcal{L}$ un fibr\'e en droites tr\`es ample et totalement sy\-m\'etri\-que sur $A$. Soit $\alpha_v$ un morphisme de Frobenius sur $A$ pour la place finie $v$. Alors, $\alpha_v$ est une isog\'enie admissible pour $\mathcal{L}$ de poids $\textnormal{q}(\alpha_v)$.
\end{prop}
\noindent \textbf{D\'emonstration} C'est la proposition 3.3. de \cite{davidhindry}.\hfill$\Box$

\section{Preuve du th\'eor\`eme \ref{semi}}

\noindent Dans la suite de ce paragraphe, on va prouver le th\'eor\`eme \ref{semi}. Pour cela nous suivrons la preuve de \cite{davidhindry} avec les modifications indiqu\'ees dans l'introduction et notamment, d'une part en introduisant dans notre cadre les id\'ees de \cite{ad2003} et d'autre part en utilisant syst\'ematiquement l'indice d'obstruction. Ce paragraphe reposant pour une forte part sur les techniques developp\'ees dans \cite{davidhindry} nous avons choisi d'insister sur les aspects qui diff\`erent mais par contre de renvoyer \`a cette r\'ef\'erence lorsqu'il n'y a pas de modification autre que typographique.

\medskip

\noindent Commen\c{c}ons tout d'abord par quelques r\'eductions.

\subsection{R\'eductions}

\noindent Quitte \`a faire une extension de degr\'e born\'e de $K$ et quitte \`a prendre une vari\'et\'e ab\'elienne isog\`ene \`a la vari\'et\'e de d\'epart, on supposera d\'esormais toujours que les hypoth\`eses de la proposition \ref{shimura} sont satisfaites. De plus une vari\'et\'e ab\'elienne de type C.M. ayant potentiellement partout bonne r\'eduction, nous supposerons \'egalement avoir choisi $K$ de sorte que $A/K$ a partout bonne r\'eduction. Ceci nous permettra notamment d'appliquer le lemme \ref{rami}.
On note 
\[\dd_n(\cdot):=\dd_{\L,K_n}(\cdot)\ \ \text{ et }\ \  \dds_n(\cdot):=[K_n:K]\dd_n(\cdot).\]
 Par ailleurs, comme on travaille avec une vari\'et\'e ab\'elienne $A/K$ de type C.M., on sait que, quitte \`a faire au d\'epart une extension born\'ee de $K$ (ce que l'on fait) de sorte que les endomorphismes de $A$ soient tous d\'efinis sur $K$, l'extension $K_n/K$ est ab\'elienne (\textit{cf.} par exemple \cite{baksi} Th 9.2 ou \cite{waterhouse} Cor. 2 du theorem 5). Enfin, quitte \`a prendre un multiple de $n$, on peut toujours supposer ce $n$ assez grand (devant les diff\'erentes constantes $c_i:=c_i(A/K,\mathcal{L})$ intervenant dans ce qui pr\'ec\`ede et dans la suite).

\subsection{L'hypoth\`ese (H)}

\noindent Soit $C_0$ une constante (ne d\'ependant que de $A/K$ et de $\mathcal{L}$ : comme dans \cite{davidhindry} $C_0$ sera prise assez grande de sorte \`a pouvoir appliquer les estimations asymptotiques (via Chebotarev) sur les ensembles $\P_r$). Soit $P$ un point de $A(\overline{K})$. On fixe deux entiers, $\rho_{\min}$ et $\rho_{\max}$, ne d\'ependant que de $g$, que nous expliciterons plus tard mais dont la valeur est fix\'ee une fois pour toute. On pose $\dds_n=\dds_n(P)$. On se donne aussi un point $Q\in A(\overline{K})$ qui v\'erifie l'\textbf{hypoth\`ese (H)} suivante :
\begin{enumerate}
\item Il existe un entier $\rho\in\{\rho_{\text{min}},\ldots,\rho_{\text{max}}\}$ tel que
\[ \h_{\L}(Q)\leq \frac{c}{\dd_n(Q)}\left(\frac{\log\log \dds_n}{C_0\log\dds_n} \right)^{(g+1)!\rho}.\]
\item Il existe une constante $c_g$ telle que $(g+1)!(2g.g!)^g\rho_{\min}\geq c_g>0$, telle que $Q=f(\sigma(P))$ o\`u $f$ est une isog\'enie de poids au plus $(C_0\log\dds_n)^{c_g}$ et $\sigma$ est un \'el\'ement de $\Gal(\overline{K}/K)$.
\end{enumerate}

\medskip

\begin{lemme}\label{deltareduit} Avec les notations pr\'ec\'edentes, si $Q$ v\'erifie l'hypoth\`ese (H), alors on a
\[\log \dd_n(Q)\leq C_0\log\dds_n(P).\]
\end{lemme}
\demo Soit $X$ une $K_n$-vari\'et\'e de dimension $d<g$ passant par le point $P$ telle que $\dd_n(P)=(\deg_{\L_{K_n}}X)^{\frac{1}{g-d}}$. La vari\'et\'e $f(\sigma(X))$ est d\'efinie sur $K_n$ et passe par $Q$, donc
\begin{align}
\dd_n(Q)	&\leq (\deg_{\L_{K_n}}f(X))^{\frac{1}{g-d}}\\
		& \leq \textnormal{q}(f)^{\frac{d}{g-d}}\dd_n(P)\label{ipl}\hspace{1cm}\text{par le lemme \ref{degre} point 2.}\\
		& \leq \textnormal{q}(f)^g\dds_n(P)
\end{align}
\noindent Or le poids $\textnormal{q}(f)$ est par l'hypoth\`ese (H) polynomial en $C_0\log\dds_n(P)$. On peut donc conclure en passant au $\log$.\hfill$\Box$

\subsection{\`A propos des param\`etres\label{para}}
\noindent On note $Q\in A(\overline{K})$ un point v\'erifiant l'hypoth\`ese (H) et $V$ une sous-$K_n$-vari\'et\'e irr\'eductible de $A_{K_n}$ de dimension minimale, r\'ealisant $\dd_n(Q)$. En utilisant un lemme de Siegel fin, le th\'eor\`eme de Bombieri-Vaaler, nous allons construire une fonction $F$, polyn\^ome homog\`ene \`a coefficients dans $\O_{K_n}$ de degr\'e $L$ en les fonctions ab\'eliennes de $A\times A$, nul \`a un ordre $\geq T_0+1$ sur l'image $i(V)$ de $V$ dans $A\times A$, le long de l'espace tangent \`a l'origine de la sous-vari\'et\'e ab\'elienne $B$ de $A\times A$ d\'efinie comme \'etant l'image de $A$ par l'application $i : x\mapsto (x,Nx)$, $N$ \'etant un param\`etre. L'entier $L$ est le degr\'e de la fonction auxiliaire, $T_0$ est l'ordre d'annulation au point $Q$ que l'on rentre dans la machine et \`a partir duquel on va extrapoler, $N$ est un param\`etre compris entre $\sqrt{L}$ et $\sqrt{2L}$. Pour des raisons techniques (\textit{cf.} le paragraphe sur l'extrapolation de \cite{davidhindry}) on choisit pour $N$ une puissance de $2$. Au vu du r\'esultat que l'on cherche \`a prouver, on fera intervenir dans le choix des param\`etres, des facteurs $\dd_n(Q)$ et des facteurs polynomiaux en $\log \dds_n$ ou $\log\log \dds_n$.

\medskip

\noindent D'autres param\`etres $N_i$ correspondant aux degr\'e des isog\'enies $\alpha_i$ de Frobenius vont intervenir. On choisit les param\`etres de sorte que
\begin{equation}\label{nul}
\h_{\L}(N\alpha_1\circ\ldots\circ\alpha_g(Q))\leq c_{3}.
\end{equation}
\noindent Enfin, on va extrapoler $g$ fois, utilisant \`a chaque fois l'extrapolation pr\'ec\'edente. On obtient ainsi des ordres d'annulation $T_i$ chacun \'etant n\'ecessairement plus petit que le pr\'ec\'edent.

\medskip

\noindent Par ailleurs, tous ces param\`etres sont choisis polynomiaux en $\dd_n(Q)$ et en $\log\dds_n$. De plus, comme $L$ donc $N^2$ va \^etre de la forme $\dd_n(Q)(\log\dds_n)^{*}$ et que on veut comme r\'esultat une minoration du type $1/\dd_n(Q)(\log \dds_n)^{*'}$, on sait par avance grace \`a l'in\'egalit\'e (\ref{nul}) que les $N_i$ seront choisis polynomiaux en $\log\dds_n$ ou plus petits.

\subsection{Lemme de Siegel}
\defi Soit $N$ un entier. Si $S$ est un sous-$\overline{\mathbb{Q}}$-espace vectoriel de $\overline{\mathbb{Q}}^{N+1}$, on d\'efinit la \textit{hauteur logarithmique de Schmidt} \cite{schmidt} de $S,$  $h_2(S)$ comme suit : sur $\mathbb{G}_m^n(\overline{\mathbb{Q}})$ on d\'efinit 
\[h_2(x_1,\ldots,x_n)=\frac{1}{d}\left(\sum_{v\in M_K^0} d_v \log\max_{1\leq
  i\leq n}\vert x_i\vert_v+ \sum_{v\in M_K^{\infty}} d_v \log \sqrt{\sum_{1\leq
  i\leq n}\vert x_i\vert^2_v}\right),\] 
\noindent o\`u $d$ et $d_v$ sont respectivement le degr\'e et les degr\'e locaux de l'extension de $\mathbb{Q}$ engendr\'ee par les coordonn\'ees $x_1,\ldots,x_n$ ; et o\`u $M_K^{\infty}$ et $M_K^0$ d\'enote respectivement l'ensemble des places infinies et des places finies de $K/\Q$ ; et o\`u on choisit comme normalisation  de $\vert.\vert_v$ pour $v\in M_K^0$ divisant $p$ la suivante : $|p|_v=p^{-1}$. Soit $N$ un entier. On d\'efinit alors la hauteur $h_2$ d'un sous-$\overline{\mathbb{Q}}$-espace vectoriel $S$ alg\'ebrique de dimension $s$ de $\overline{\mathbb{Q}}^{N+1}$ par :
\[h_2(S)=h_2(\mathbf{x}_1\wedge\ldots\wedge \mathbf{x}_s),\]
\noindent o\`u $\mathbf{x}_1,\ldots, \mathbf{x}_s$ est une base de $S$ sur un corps de nombres quelconque sur lequel $S$ est d\'efini.

\medskip

\noindent Rappelons le lemme de Siegel que nous allons utiliser : il s'agit du th\'eor\`eme de Bombieri-Vaaler \cite{B-V}. Il y a ici une diff\'erence fondamentale avec l'article de David et Hindry : on veut pouvoir controler la d\'ependance en l'extension $K_n/K$. 

\medskip

\begin{theor}\label{bombierivaaler}\textnormal{\textbf{(Bombieri-Vaaler)}} Soit $\mathbb{F}/\Q$ un corps de nombres de degr\'e $d$. Soient $M$ et $N$ deux entiers strictement positifs et $S$ un sous-$\mathbb{F}$-espace vectoriel de $\mathbb{F}^N$ de dimension $N-M>0$. Il existe un \'el\'ement non nul $\textbf{x}$ dans $\O_\mathbb{F}^N$ de $S$ tel que 
\[h_2(1,\textbf{x})\leq \frac{1}{2d}\log\mid\text{Disc}(\mathbb{F}/\Q)\mid+\frac{1}{N-M}h_2(S).\]
\end{theor}

\medskip

\noindent Avant de poursuivre, faisons quelque brefs rappels concernant la notion de vari\'et\'e projectivement normale.

\medskip

\defi On dit qu'une sous-vari\'et\'e $X$ de $\mathbb{P}_n$ est \textit{projectivement normale} si son anneau de coordonn\'ees $S(X)$ est un anneau normal (\textit{i.e.}, int\'egralement clos).

\medskip

\noindent On peut montrer (\textit{cf.} par exemple Birkenhake-Lange \cite{lange} p. 190-193) que $X\subset \mathbb{P}_n$ est projectivement normale si et seulement si elle est normale, et pour tout $d\geq 0$ la fl\`eche naturelle 
\[ H^0(\mathbb{P}_n,\mathcal{O}_{\mathbb{P}_n}(d))\rightarrow H^0(X,\mathcal{O}_X(d))\]
\noindent est surjective.

\medskip

\noindent Concernant les vari\'et\'es ab\'eliennes plong\'ees de mani\`ere projectivement normale, on a le r\'esultat suivant que l'on trouve par exemple dans \cite{lange} theorem 3.1 p. 190.

\begin{prop}Soient $A/K$ une vari\'et\'e ab\'elienne, et $\mathcal{L}$ un fibr\'e en droites ample sur $A$. Pour tout $n\geq 3$, le fibr\'e $\mathcal{L}^{\otimes n}$ d\'efinit un plongement projectivement normal de $A$ dans un espace projectif $\mathbb{P}_r$.
\end{prop}

\medskip

\noindent Soient $A/K$ une vari\'et\'e ab\'elienne sur un corps de nombres, munie d'un fibr\'e en droites sym\'etrique ample $\mathcal{L}$. Quitte \`a travailler avec $\mathcal{L}^{\otimes 4}$ plut\^ot qu'avec $\mathcal{L}$, on peut supposer que $\mathcal{L}$ est tr\`es ample, totalement sym\'etrique et d\'efinit un plongement projectivement normal de $A$ dans un projectif $\mathbb{P}_n$. On note $\mathcal{M}=\mathcal{L}\boxtimes\mathcal{L}$ le fibr\'e en droites sur $A\times A$ associ\'e \`a $\mathcal{L}$. On va maintenant pouvoir construire la fonction $F$ recherch\'ee.

\medskip

\noindent Soient $L$ et $T$ deux entiers positifs. On note $\left\{s_0,\ldots,s_l\right\}$ une base de $H^0(A\times A,\mathcal{M})$. On peut, par projective normalit\'e, choisir une base $\left\{Q_1,\ldots,Q_m\right\}$ du $K$-vectoriel $H^0\left(A\times A,\mathcal{M}^{\otimes L}\right)$ telle que tous les $Q_i$ sont homog\`enes de degr\'e $L$ en les $s_j$. De plus, on peut aussi voir les $s_i$ comme des $(1,1)$-formes homog\`enes de $K[\textbf{X},\textbf{Y}]$ o\`u $\textbf{X}=(X_0,\ldots,X_n)$, et $\textbf{Y}=(Y_0,\ldots,Y_n)$. Enfin on note $T_B$ l'espace tangent \`a l'origine de la sous-vari\'et\'e ab\'elienne $B=i(A)$ de $A\times A$ d\'efinie par $y=[N]x$. 

\medskip

\noindent \textbf{But :} fabriquer un polyn\^ome, $F=\sum_{i=1}^m b_iQ_i$, \`a coefficients entiers dans $\mathcal{O}_{K_n}$, en les fonctions ab\'eliennes de $A\times A$, tel que $F$ est de ``petite'' hauteur, et tel que $F$ s'annule \`a un ordre sup\'erieur \`a $T_0+1$ sur $i(V)$, le long de $T_B$.

\medskip

\noindent En notant $\Theta$ l'application th\^eta d\'efinie sur $T_{A(\mathbb{C})}$ par la composition 
\[
\xymatrix{
T_{A(\mathbb{C})}\ \ar[r]^{\textnormal{exp}_{A(\mathbb{C})}}	& \  A(\mathbb{C}) \ar[r]^{\varphi_{\mathcal{L}}}	& \ \mathbb{P}_n(\mathbb{C}) }
\]
\noindent associ\'ee \`a $\mathcal{L}$, ceci correspond \`a trouver une solution de petite hauteur au syst\`eme d'inconnues les $b_i$
\begin{equation}\label{systeme1}
 \frac{\partial^{\kappa}F\left(\Theta(\textbf{u}+\textbf{z}),\Theta(N(\textbf{u}+\textbf{z})\right)}{\partial \textbf{z}^{\kappa}}_{\textnormal{\Large{$\mid$}}\textbf{z}=0}=0,
\end{equation}
\noindent pour tout $\mid \kappa\mid\leq T$ et $\textbf{u}\in T_{A(\mathbb{C})}$ tels que $\Theta(\textbf{u})\in V(\overline{K})$. 

\medskip

\noindent On reprend le lemme 5.1. de \cite{davidhindry} en rempla\c{c}ant $K$ par $K_n$. On obtient alors directement : 

\medskip

\begin{lemme}\label{rang}Il existe une constante strictement positive $c_{4}$ telle que le rang du syst\`eme (\ref{systeme1}) sur $K_n$ est major\'e par
\[\textbf{rg}=c_{4}\left(T_0\delta_n(Q)\right)^{g-d_0}\left(LN^2\right)^{d_0}.\]
\end{lemme}

\medskip

\noindent On peut maintenant construire la fonction que l'on veut. \'Etant donn\'e un polyn\^ome $F$ \`a coefficients dans $\overline{K}$, on note $h(F)$ la hauteur logarithmique absolue du point projectif d\'efini par $1$ et les coefficients de $F$.

\medskip

\begin{lemme}\label{sieg}Si $T_0\dd_n(Q)<2L^2$, alors il existe une fonction $F$ solution du syst\`eme (\ref{systeme1}), de hauteur major\'ee par
\[h(F)\leq c_{5}\left(\log [K_n:K]\right)^2+\frac{\textbf{rg}\times\left(C_0T_0\log(\dd_n(Q))+L\right)}{L^{2g}}.\]
\end{lemme}
\demo Il suffit de reprendre la preuve du lemme 5.4. de \cite{davidhindry} et d'appliquer le th\'eor\`eme \ref{bombierivaaler} en lieu et place du classique lemme de Siegel. Comme $K_n=K(A[n])$, on utilise le corollaire \ref{coro} pour majorer le discriminant appara\^issant dans le th\'eor\`eme \ref{bombierivaaler}. C'est ce discriminant qui nous donne le terme $ c_{5}\left(\log [K_n:K]\right)^2$.\hfill$\Box$

\subsection{Lemme de z\'eros}
\noindent Avec les notations pr\'ec\'edentes, on pose $\S=\left\{\sigma(Q)\ |\ \sigma\in \Gal(\overline{K}/K_n)\right\}$. Pour chaque premier $p$, fixons une place $v_p$ de $K_n$ au-dessus de $p$ et notons $NR_{K_n}$ l'ensemble des telles places $v_p$ telles que $v_p|p$ est non-ramifi\'ee dans $K_n$. Pour $r$ compris entre $1$ et $g$, posons de plus
\[\P_r=\left\{\text{Id}\right\}\cup\left\{\alpha_v\ |\ v\in NR_{K_n},\ \ \frac{N_r}{2}\leq N(v)\leq N_r\right\},\]
\noindent o\`u les $N_r$ sont des param\`etres qui seront sp\'ecifi\'es ulterieurement. Il faut notamment les choisir de sorte que $\frac{N_i}{\log N_i}\geq 2\log\dds_n \geq \frac{4}{\log 2}\log [K_n:K]$ afin d'\^etre sur que les ensembles $\P_r$ soient non vides (\textit{cf.} les points (2) des lemmes \ref{nombrepremier} et \ref{rami}). En supposant ceci v\'erifi\'e, le th\'eor\`eme de Chebotarev nous assure alors que 
\[\textnormal{Card}(\P_r)\geq\frac{c_{6}N_r}{\log N_r}.\]
\noindent \textbf{Notation.} Si $\alpha$ est un \'el\'ement de $\mathcal{P}_r$ pour un certain entier positif $r$, associ\'e au nombre premier $p$, on notera $\tilde{\alpha}$ l'op\'erateur $\alpha\circ \Phi_p^{-1}$, compos\'e d'une isog\'enie et d'un morphisme de $\Gal(K_n/K)$. Par convention, si $\alpha$ est l'identit\'e, on prendra \'egalement l'identit\'e pour $\widetilde{\alpha}$.

\medskip

\noindent Pour tout entier $r$ entre $1$ et $g$, on note alors
\[\S^{(r)}=\left\{\aa_r\circ\cdots\circ\aa_{g}(Q))\ | \ Q\in \S,\ \alpha_i\in \P_i\right\}.\]
\noindent On convient que $\S^{(g+1)}=\S$.

\medskip

\noindent On rappelle une variante (affaiblie) du lemme de z\'eros d\'emontr\'e dans \cite{davidhindry} (th\'eor\`eme 4.1.).

\medskip

\begin{theo}\label{zero}\textnormal{\textbf{(Lemme de z\'eros)}} On utilise les notations pr\'ec\'edemment introduites. Soient $A/K$ une vari\'et\'e ab\'elienne de dimension $g$, plong\'ee dans un espace projectif $\mathbb{P}_m$ de fa\c{c}on projectivement normale et $M$ un entier strictement positif. On se donne une forme $F\in K_n[X_0,\ldots,X_m]$ homog\`ene de degr\'e $L$, non-identiquement nulle sur $A_{K_n}$. On suppose que $F$ s'annule \`a un ordre sup\'erieur \`a $1+M$ le long de $T_{A_{K_n}}$ en tous les points de $\S^{(1)}$. Sous ces hypoth\`eses, 

il existe un entier $r\in\{1,\ldots,g\}$, une sous-$K_n$-vari\'et\'e de $A_{K_n}$, $V$, stricte et $K_n$-irr\'e\-duc\-tible, de dimension $d\geq g-r$, telle que $V_{\overline{K}}$ contient un \'el\'ement de $\S^{(r+1)}$, incompl\`e\-tement d\'efinie dans $A$ avec multiplicit\'e sup\'erieure \`a $\frac{1}{g}M$ le long de l'espace tangent \`a l'origine $T_{A_{K_n}}$ par des formes de degr\'e inf\'erieur \`a $2LN_1\times\cdots\times N_{r-1}$, telle que :
\[M^{g-d}\deg_{\L_{K_n}}\left(\bigcup_{\alpha\in \P_r}\widetilde{\alpha}(V)\right)\leq c_{7}(LN_1\times\cdots\times N_{r-1})^{g-d}.\]
\end{theo}
\demo C'est le th\'eor\`eme 4.1. (version galoisienne) et la scolie 4.8. de \cite{davidhindry} appliqu\'es avec 
\[K=K_n,\hspace{1cm} \forall i\in\{1,\ldots,g\}\ \ T_i=\frac{1}{g}M,\hspace{1cm}  L=L,\hspace{1cm} \mathcal{V}=T_{A_{K_n}}.\]
\noindent Plus pr\'ecis\'ement, on reprend leur preuve,la seule diff\'erence \'etant la suivante : l\`a o\`u ils introduisent (p.27, paragraphe 4.2.) la suite d'id\'eaux
\[\mathfrak{I}_1=(P),\ \ \ \mathfrak{I}_{r+1}=(\partial^{T_r}_{0,\alpha}\mathfrak{I}_r ;\ \alpha\in\mathcal{P}_r),\]
\noindent nous introduisons la suite
\[\mathfrak{I}_1=(P),\ \ \ \mathfrak{I}_{r+1}=(\partial^{T_r}_{0,\alpha_p}\Phi_p(\mathfrak{I}_r) ;\ \widetilde{\alpha}_p=\alpha_p\circ\Phi_p^{-1}\in\mathcal{P}_r).\]
\noindent La preuve est alors la m\^eme. \hfill$\Box$

\medskip

\noindent Comme dans \cite{davidhindry}, nous aurons \'egalement besoin d'un r\'esultat suppl\'ementaire, raffinant l'in\'egalit\'e de B\'ezout :

\medskip

\begin{lemme}\label{bezout}Il existe une constante $c_8$ telle que la propri\'et\'e suivante est vraie : soient $V$ une sous-$K_n$-vari\'et\'e stricte de $A$, irr\'eductible sur $K_n$ et $F$ une forme de degr\'e $\mu$ sur $A$, d\'efinie sur $K_n$. On suppose que $F$ est nulle avec multiplicit\'e sup\'erieure \`a $m$ sur une sous-$K_n$-vari\'et\'e $V'$ de $A$, et on suppose \'egalement que les intersections $V\cdot V'$ et $V\cdot\mathcal{Z}(F)$ ont une composante W, $K_n$-irr\'eductible en commun, de codimension $1$ dans $V$. Dans ce cas, on a l'in\'egalit\'e
\[\deg_{\L_{K_n}}W\leq c_8\frac{(\deg_{\L_{K_n}} V)\mu}{m}.\]
\end{lemme}
\demo Il s'agit du lemme 4.9. de \cite{davidhindry} avec $K=K_n$.\hfill$\Box$

\subsection{Extrapolation}
\noindent L'extrapolation suit pour partie ce qui est fait par David et Hindry dans \cite{davidhindry}. Ceci dit, il y a tout de m\^eme une diff\'erence : comme la fonction auxiliaire n'est plus construite \`a coefficients entiers mais \`a coefficients dans $\O_{K_n}$, on va \`a la mani\`ere des cas non-ramifi\'e de \cite{AZ} et \cite{art3}, faire intervenir les automorphismes de  Frobenius $\Phi_p$ de l'extension ab\'elienne $K_n/K$. Ceci \'etant not\'e, les choses fonctionnent bien ensuite \`a condition de se restreindre au cas des premiers non-ramifi\'es dans $K_n$. C'est pour cela que l'on d\'efinit les ensembles $\P_r$ tels qu'on les d\'efinit. On \'enonce le lemme qui nous permet d'extrapoler. Notons que l'on utilise ici le fait que l'extension $K_n/K$ est ab\'elienne.

\medskip

\begin{lemme}\label{nonramifie} Soient $x\in \mathcal{O}_{K_n}$, $p$ un nombre premier non-ramifi\'e dans $K_n$ et $v$ une valuation sur $\overline{K}$ \'etendant $p$. En notant $\Phi_p\in\textnormal{Gal}(K_n/K)$ l'automorphisme de Frobenius associ\'e \`a $p$, on a 
\[\mid x^p-\Phi_p x\mid_v\leq p^{-1}.\]
\end{lemme}
\demo C'est le lemme 3.1. de \cite{AZ}.\hfill$\Box$

\medskip

\rem Notons que c'est cette restriction, le fait de se restreindre \`a ne travailler qu'avec les premiers non-ramifi\'es dans $K_n$, qui fait appara\^itre les facteurs $\log\dds_n$ plut\^ot que $\log\dd_n(P)$. Pour corriger ceci, on pourrait s'inspirer de l'article \cite{AZ} (et \cite{art3}), mais, ind\'ependamment des autres complications, il y aurait cette fois-ci non plus $g$, mais $2^g$ \'etapes d'extrapolation \`a faire.

\medskip

\noindent On reprend le paragraphe 6 de \cite{davidhindry} et on remplace dans le th\'eor\`eme 6.4, le corps $K(Q)$ par $K_n(Q)$. Ceci se fait sans autre changement et on obtient ainsi~:

\medskip

\begin{prop}\label{extrapol} Pour tout $i$ compris entre $1$ et $g$, on fait les hypoth\`eses suivantes sur les param\`etres :
\[N^2\left(\prod_{i=1}^gN_i\right)\h_{\L}(Q)\leq c_{3},\]
\noindent et 
\[T_i\log N_{g-i}\geq 2L\text{ et si }i\leq g-1,\ \  T_i\log N_{g-i}> 2c_{9}\left(T_{i+1}\log(T_{i+1}+L)+L+h(F)\right).\]
\noindent Dans ce cas, la fonction auxiliaire $F$ est nulle le long de $T_{B(\C)}$ \`a un ordre sup\'erieur \`a $T_g$ en tout point de $\S^{(1)}$.
\end{prop}
\demo Il s'agit de reprendre les calculs de la proposition 6.5. de \cite{davidhindry} (elle m\^eme bas\'ee sur l'extrapolation de Laurent \cite{laurent}) en passant de $K$ \`a $K_n$. On conserve donc leurs notations. Il s'agit de la m\^eme chose que ce qui est fait dans les propositions 6.1 et 7.1 de \cite{art3} dans le cadre des courbes elliptiques. On raisonne par r\'ecurrence, on suppose $F$ nulle \`a un ordre sup\'erieur \`a $T_i$ le long de $T_{B(\C)}$ en tous les points de $\S^{(g+1-i)}$, et on va montrer la m\^eme chose \`a au rang $i+1$. Plus exactement, \'etant donn\'e un point $R_{g+1-i}$ de $\S^{(g+1-i)}$, plut\^ot que de montrer que $F$ s'annule en un point $\aa_v(R_{g+1-i})$ \`a un ordre sup\'erieur \`a $T_{i+1}$ avec $\alpha_v\in\P_{g-i}$, on va montrer, ce qui est \'equivalent, que $\Phi_p(F)$ s'annule au point $\a_v(R_{g+1-i})$ au m\^eme ordre. Cette reformulation nous permettra d'appliquer le lemme \ref{nonramifie}. Notons que le rang $i=0$, autrement dit l'initialisation de la r\'ecurrence, est vrai par construction de $F$ par le lemme de Siegel. 

\medskip

\noindent Soit $v$ la place de $K_n$ dans $NR_{K_n}$ correspondant \`a $\alpha_v\in\P_{g-i}$. Elle est d'uniformisante $\pi_v$. Soit $R$ un point de $\S^{(g+1-i)}$ d\'efini sur une extension $K_n'/K_n$, et soit $w$ une place de $K_n'$ au dessus de $v$. Notons $\textbf{R}=(R_0,\ldots,R_n)$ un syst\`eme de coordonn\'ees projectives de $R$ dans $\mathcal{O}_w$, telles que $\mid\mid \textbf{R}\mid\mid_w=1$. Soit $\partial^{\kappa}$ un op\'erateur diff\'erentiel d'ordre $\mid\kappa\mid\leq T_{i+1}$ le long de $T_{B(\mathbb{C})}$. L'application de l'hypoth\`ese de r\'ecurrence (annulation de $F$ \`a un ordre sup\'erieur \`a $T_i$ le long de $T_{B(\C)}$ en tous les points de $\S^{(g+1-i)}$), et l'application du petit th\'eor\`eme de Fermat dans le cadre des vari\'et\'es ab\'eliennes conduisant \`a l'in\'egalit\'e (20) p.47 de \cite{davidhindry} et le lemme \ref{nonramifie} nous donnent 
\begin{equation}\label{vingt}
 \left | \Phi_p\left(\partial^{\kappa}F\right)\left(\mathbf{F}_{\alpha_v}(\mathbf{R}),\mathbf{F}^{(N)}\circ\mathbf{F}_{\alpha_v}(\mathbf{R})\right)\right |_w\leq |\pi_v|_w^{T_i-\mid\kappa\mid},
\end{equation}
\noindent o\`u $\mathbf{F}_{\alpha_v}$ et $\mathbf{F}^{(N)}$ sont des formes homog\`enes de $\mathcal{O}_K[\mathbf{X}]$ de degr\'e respectifs $\textnormal{N}(v)$ et $4^{m+1}$, repr\'esentant respectivement l'endomorphisme de Frobenius sur $A$ associ\'e \`a $v$, et la multiplication par $N=2^{m+1}$.

\medskip

\noindent On veut maintenant sommer sur toutes les places $w$ au-dessus de $v$. Malheureusement, le choix du syst\`eme de coordonn\'ees projectives pour $R$ d\'e\-pend de $w$. On est donc oblig\'e d'alourdir les notations pour pallier ce probl\`eme. Soient $S$, $S_N$, $S_{\alpha_v}$, $S_{N,\alpha_v}$ des coordonn\'ees projectives non nulles de $R$, $\mathbf{F}^{(N)}(\mathbf{R})$, $\mathbf{F}_{\alpha_v}(\mathbf{R})$, $\mathbf{F}^{(N)}\circ\mathbf{F}_{\alpha_v}(\mathbf{R})$ respectivement. On note de plus $S_{w,N}$, $S_{w,\alpha_v}$, $S_{w,N,\alpha_v}$ des coordonn\'ees de ces points de valeur absolue $w$-adique maximale.

\medskip

\noindent Soit maintenant $\partial^{\kappa}$ un op\'erateur diff\'erentiel de longueur minimale pour lequel
\[ \Phi_p\left(\partial^{\kappa}F\right)\left(\mathbf{F}_{\alpha_v}(\mathbf{R}),\mathbf{F}^{(N)}\circ\mathbf{F}_{\alpha_v}(\mathbf{R})\right)\]
\noindent est non nul. Si $|\kappa|$ est sup\'erieur \`a $T_{i+1}$, on a gagn\'e. Sinon on applique la formule de Leibniz en utilisant que $F$ est bihomog\`ene de bidegr\'e $(L,L)$. On a donc 
\begin{align*}
\Phi_p\left(\partial^{\kappa}F\right)\left(\frac{\mathbf{F}_{\alpha_v}(\mathbf{R})}{S_{\alpha_v}},\frac{\mathbf{F}^{(N)}\left(\mathbf{F}_{\alpha_v}(\mathbf{R})\right)}{S_{N,\alpha_v}}\right) & =  \frac{\Phi_p\left(\partial^{\kappa}F\right)\left(\mathbf{F}_{\alpha_v}(\mathbf{R}),\mathbf{F}^{(N)}\left(\mathbf{F}_{\alpha_v}(\mathbf{R})\right)\right)}{S_{\alpha_v}^LS_{N,\alpha_v}^L}
\end{align*}	
\noindent Or ceci est \'egal \`a 
\[ \Phi_p\left(\partial^{\kappa}F\right)\left(\frac{\mathbf{F}_{\alpha_v}(\mathbf{R})}{S_{w,\alpha_v}},\frac{\mathbf{F}^{(N)}\left(\mathbf{F}_{\alpha_v}(\mathbf{R})\right)}{S_{w,N,\alpha_v}}\right)\cdot \frac{\left(S_{w,\alpha_v}S_{w,N,\alpha_v}\right)^L}{\left(S_{\alpha_v}S_{N,\alpha_v}\right)^L}.\]
\noindent On r\'e\'ecrit alors l'in\'egalit\'e $(\ref{vingt})$ en passant au log, en sommant sur toutes les places $w$ au-dessus de $v$ et en notant $d_w$ les degr\'es locaux :
\begin{eqnarray}\label{e12}
\lefteqn{\sum_{w|v}d_w\log\left(\left|\Phi_p\left(\partial^{\kappa}F\right)\left(\frac{\mathbf{F}_{\alpha_v}(\mathbf{R})}{S_{\alpha_v}},\frac{\mathbf{F}^{(N)}\left(\mathbf{F}_{\alpha_v}(\mathbf{R})\right)}{S_{N,\alpha_v}}\right)\right|_w\right)}\\
& & \leq \left(T_i-|\kappa|\right)\sum_{w|v}d_w\log(|\pi_v|_w)+L\sum_{w|v}d_w\log\left(\frac{|S_{w,\alpha_v}S_{w,N,\alpha_v}|_w}{|S_{\alpha_v}S_{N,\alpha_v}|_w}\right)\label{eqnar}.
\end{eqnarray}
\noindent Or 
\begin{equation*} \label{sys21}
\sum_{w|v}d_w\log(|\pi_v|_w)=[K_n':K_n]\log(|\pi_v|_v)= -\frac{[K_n':K_n]}{d_v}\log\left(\textnormal{N}(v)\right)
\end{equation*}
\noindent De plus, on peut voir que 
\begin{align*}
\sum_{w|v}d_w\log\left(\frac{|S_{w,\alpha_v}S_{w,N,\alpha_v}|_w}{|S_{\alpha_v}S_{N,\alpha_v}|_w}\right)& \leq [K_n':K_n]\left(h_{\mathcal{L}}(\alpha_v(R))+h_{\mathcal{L}}(N\alpha_v(R))\right)\\
	& \leq [K_n':K_n]\left(\textnormal{N}(v)\widehat{h}_{\mathcal{L}}(R)\!+N^2\textnormal{N}(v)\widehat{h}_{\mathcal{L}}(R)\!+c_{14}\right).
\end{align*}
\noindent C'est l'in\'egalit\'e $(21)$ p. 49 de \cite{davidhindry}. On obtient ainsi, en notant $\text{G}$ le membre de gauche de l'in\'egalit\'e (\ref{e12}) la majoration suivante :
\begin{align}\label{tnl} 
\text{G}	& \leq -c_{10}T_i[K_n':K_n]\log N_{g-i}+L[K_n':K_n](N(v)\left(N^2+1)\h_{\L}(R)+c_{14}'\right)\\
		& \leq  -c_{10}T_i[K_n':K_n]\log N_{g-i}+c_{11}L[K_n':K_n]\leq  -c_{12}T_i[K_n':K_n]\log N_{g-i},
\end{align}
\noindent par choix des param\`etres et par l'hypoth\`ese. Pr\'ecis\'ement c'est pour obtenir cette derni\`ere in\'egalit\'e que l'on applique l'hypoth\`ese \textbf{(H)}, via l'in\'egalit\'e (\ref{nul}) suppos\'ee v\'erifi\'ee dans les hypoth\`eses de la proposition. Il reste \`a majorer le terme $-\text{G}$. Par d\'efinition de la hauteur (absolue logarithmique) projective, on a 
\begin{align*}
\frac{-1}{[K_n':K_n]}\text{G} & \leq h\left(\left(\Phi_p\left(\partial^{\kappa}F\right)\left(\frac{\mathbf{F}_{\alpha_v}(\mathbf{R})}{S_{\alpha_v}},\frac{\mathbf{F}^{(N)}\left(\mathbf{F}_{\alpha_v}(\mathbf{R})\right)}{S_{N,\alpha_v}}\right)\right)^{-1}\right)\\
  &  = h\left(\Phi_p\left(\partial^{\kappa}F\right)\left(\frac{\mathbf{F}_{\alpha_v}(\mathbf{R})}{S_{\alpha_v}},\frac{\mathbf{F}^{(N)}\left(\mathbf{F}_{\alpha_v}(\mathbf{R})\right)}{S_{N,\alpha_v}}\right)\right),
\end{align*}
\noindent ceci ayant un sens grace \`a l'hypoth\`ese de non nullit\'e de $\Phi_p\left(\partial^{\kappa}F\right)(\cdots)$. Il ne reste maintenant plus qu'\`a majorer cette derni\`ere hauteur. Il s'agit d'un calcul classique (\textit{cf.} par exemple \cite{davidhindry} p. 50). On obtient
\begin{equation}\label{gauche}
\frac{-1}{[K_n':K_n]}\text{G}\leq c_{13}\left(T_{i+1}\log(T_{i+1}+L)+LN^2\textnormal{N}(v)\widehat{h}_{\L}(R)+h(F)\right).
\end{equation}
\noindent Finalement, en mettant ensemble les in\'egalit\'es $(\ref{tnl})$ et $(\ref{gauche})$, on obtient
\begin{equation}\label{extra1}
T_i\log N_{g-i}\leq  c_{14}\left(T_{i+1}\log(T_{i+1}+L)+L+h(F)\right).
\end{equation}

\noindent Les hypoth\`eses permettent alors de conclure.\hfill $\Box$

\subsection{Choix complet des param\`etres}
Jusque l\`a on a du imposer les in\'egalit\'es suivantes sur les param\`etres :

\[\forall i, N_i\gg(\log\dds_n)^2,\hspace{1cm} N^2\sim L,\hspace{.5cm}\text{ et si }\ \ i\leq g-1,\  T_0\dd_n(Q)\ll L^2\ll (T_i\log\log\dds_n)^2.\]

De plus, $L$ \'etant lin\'eaire en $\dd_n(Q)$ (\textit{cf.} le paragraphe \ref{para}), il en est de m\^eme pour les $T_i$ qui sont strictement d\'ecroissants, en fait v\'erifiant :
\[T_{i+1}\leq T_i\frac{\log\log \dds_n}{C_0\log \dds_n}.\]
\noindent On doit \'egalement avoir $T_0\geq\ldots\geq T_{g-1}\geq \frac{(\log\dds_n)^2}{\log\log\dds_n}$. Ainsi, on a :

\[T_g\simeq T_0\left(\frac{\log\log\dds_n}{C_0\log\dds_n}\right)^{g}=\dd_n(Q)\frac{T_0}{\dd_n(Q)}\left(\frac{\log\log\dds_n}{C_0\log\dds_n}\right)^{g}.	\]
\noindent Enfin la derni\`ere chose \`a v\'erifier est la suivante : 
\[\forall i\in\{1,\ldots,g-1\},\ \ C_0T_0^{g+1-d_0}\dd_n(Q)^{g-d_0}\log\dds_n \ll cT_iL^{2(g-d_0)}\log\log \dds_n.\]
\noindent On veut (afin d'obtenir un r\'esultat optimal vis \`a vis de la m\'ethode) que $L$ soit le plus petit possible, ceci nous permet d'en d\'eduire les valeurs optimales pour les $T_i$ et pour $L$ : 
\[ L=\left[C_0^{2g}\dd_n(Q)(\log\dds_n)^{2g-1}(\log\log\dds_n)^{-2g}\right], N=2^{m+1}, m=\left[\frac{\log L}{2\log 2}\right],\]
\noindent et, pour tout $i\in \{1,\ldots,g\}$
\[T_0=\left[C_0^{3g-\frac{1}{2}}\dd_n(Q)(\log\dds_n)^{3g-2}(\log\log\dds_n)^{-3g}\right], T_i=\left[T_0\left(\frac{\log\log\dds_n}{C_0\log\dds_n}\right)^{i}\right].\]
\noindent Par ailleurs on pose

\[\forall i\in\{1,\ldots,g\}\hspace{.5cm} N_i=\left(\frac{C_0\log\dds_n}{\log\log\dds_n}\right)^{i.i!\rho},\]
\noindent pour un certain param\`etre $\rho$ entre $\rho_{\text{min}}$ et $\rho_{\text{max}}$, avec $\rho_{max}=(g+2)(2g.g!)^g\rho_{\min}$ et $\rho_{\min}=2g+5$. Avec ces choix de param\`etres on v\'erifie que l'in\'egalit\'e (\ref{nul}) est bien v\'erifi\'ee~: 
\[\h_{\L}(N\alpha_1\circ\ldots\circ\alpha_g(Q))\leq c_{3}.\]

\subsection{Fin de l'extrapolation}\label{fin}
\noindent On consid\`ere l'application 

\[
\begin{array}{ccccccccc}
\phi : & A & \overset{i}{\hookrightarrow} & A\times A  & \hookrightarrow & \mathbb{P}_n\times\mathbb{P}_n & \underset{\textnormal{Segre}}{\hookrightarrow} & \mathbb{P}_{(n+1)^2-1}\\
	& x & \mapsto         & (x,[N]x) &                 &          &                                          &
\end{array}
\]

\begin{prop}\label{resume}Soient $\rho$ compris entre $\rho_{\text{min}}$ et $\rho_{\text{max}}$ et $Q$ un point de $A(\overline{K})$ v\'erifiant l'hypoth\`ese \textnormal{(H)} pour $\rho$. Il existe une $K_n$-vari\'et\'e $V$, $K_n$-irr\'eductible, stricte de $A$, de dimension $d$, telle que : il existe un point $Q_1\in A(\overline{K})$, de la forme $Q_1=\aa_{r+1}\circ\ldots\circ\aa_g(Q)$, pour un certain $r\in\{1,\ldots,g\}$ et certains $\alpha_i\in\P_i$, $i\in\{r+1,\ldots,g\}$, tel que $V$ est de dimension $d\geq g-r$, que $Q_1\in V(\overline{K})$, et si $V$ n'est pas une sous-vari\'et\'e torsion, on a
\[\deg_{\L_{K_n}}V	\leq c_{15}\frac{C_0\log\log\dds_n}{N_r}\left(\frac{L^2N_1\times\cdots\times N_{r-1}}{T_g}\right)^{g-d}.\]
\noindent De plus $V$ est incompl\`etement d\'efinie par des formes de degr\'e au plus
\[c_{16}LN^2N_1\times\cdots\times N_{r-1}\]
\noindent avec multiplicit\'e sup\'erieure \`a $\frac{1}{g}T_g$.
\end{prop}
\demo On applique le th\'eor\`eme \ref{zero} \`a la fonction auxiliaire $F$ construite pr\'ec\'e\-demment, tir\'ee en arri\`ere par $\phi$. Comme $F$ est une forme bi-homog\`ene de bidegr\'e $(L,L)$ non-iden\-ti\-que\-ment nulle sur $A_{K_n}\times A_{K_n}$, elle n'est pas identiquement nulle sur $B_{K_n}$ par choix du param\`etre $N$ (on a pris $N^2\geq L+1$). De plus la proposition \ref{extrapol} nous indique que $F$ est nulle le long de $T_{B(\C)}$ \`a un ordre sup\'erieur \`a $T_g$ en tous les points de $\S^{(1)}$. La forme $G=F\circ\phi$, qui est une forme de degr\'e $L(N^2+1)$ v\'erifie donc les hypoth\`eses du lemme de z\'eros. Il suffit de v\'erifier l'in\'egalit\'e annonc\'ee pour le degr\'e de $V$. On suppose donc que $V$ n'est pas de torsion.

\medskip

\noindent On note $\tilde{V}$ une composante irr\'eductible de $V_{\overline{K}}$. Comme $V$ n'est pas une sous-vari\'et\'e de torsion, le point 1. du lemme \ref{distinct} nous dit alors qu'il n'existe pas de triplet $(\alpha,\beta,\sigma)\in \P_r^2\times\Gal(\overline{K}/K_n)$ tel que $\alpha(\tilde{V})=\beta(\sigma(\tilde{V}))$. Le choix de $N_r\geq (\log\dds_n)^{2}$ fait que la suite du calcul de la majoration du degr\'e de $V$ marche comme dans \cite{davidhindry} p.54-55, en rempla\c{c}ant $K$ par $K_n$ :
\[\deg\left(\bigcup_{\alpha\in\mathcal{P}_r,\ \sigma\in\Gal(\overline{K}/K_n)}\aa(\sigma(\tilde(V))\right)\geq \sum_{\alpha\in\mathcal{P}_r}\sum_{i=1}^M\frac{\textnormal{q}(\alpha)^d\deg(\sigma_i(\tilde{V}))}{\text{Card}\left( \text{Ker}\alpha\cap G_{\sigma_i(\tilde{V})}\right)},\]
\noindent o\`u l'on a not\'e $\sigma_i(\tilde{V})$ les diff\'erents conjugu\'ees de $\tilde{V}$. On poursuit alors comme dans \cite{davidhindry} p.55. Notamment on a
\[\mid G_{\tilde{V}}:G_{\tilde{V}}^0\mid\leq \deg \tilde{V}(4L^2N_1\times\cdots\times N_{r-1})^{g-s}\leq C_0^{\star_1}(\dds_n)^{\star_2},\]
\noindent o\`u $\star_1$ et $\star_2$ sont des constantes explicites ne d\'ependant que de $g$. De plus par choix du param\`etre $N_r$, le cardinal de l'ensemble $\P_r$ est sup\'erieur \`a $\frac{N_r}{C_0\log\log\dds_n}$. Finalement tous calculs faits, on obtient
\[\deg\left(\bigcup_{\alpha\in\mathcal{P}_r,\ \sigma\in\Gal(\overline{K}/K_n)}\aa(\sigma(\tilde(V))\right)\geq \frac{c'_{15}}{C_0}\frac{N_r\deg_{\L_{K_n}}V}{\log\log\dds_n}.\]
\noindent En rempla\c{c}ant dans l'in\'egalit\'e fournie par le lemme de z\'eros, on obtient :
\[\deg_{\L_{K_n}}V\leq c_{15}\frac{C_0\log\log\dds_n}{N_r}\left(\frac{L^2N_1\times\cdots\times N_{r-1}}{T_g}\right)^{g-d}.\]
\noindent Ceci conclut. \hfill$\Box$

\medskip

\rem Notons que en rempla\c{c}ant dans la proposition les param\`etres par leur valeur, on obtient,
\[\deg_{\L_{K_n}}V	\leq c_{15}\frac{C_0^{(2g+\frac{3}{2}+\rho(r!-1))(g-d)}\delta_n(Q)^{g-d}(\log\dds_n)^{(2g+\rho(r!-1))(g-d)}}{N_r(\log\log\dds_n)^{(2g+\rho(r!-1))(g-d)-1}}.\]

\medskip

\begin{cor}\label{cor2}Avec les notations de la proposition pr\'ec\'edente et si $V$ n'est pas de torsion, on a l'in\'egalit\'e
\[\dd_n(Q_1)\leq \frac{C_0^{2g+\frac{3}{2}}(\log\dds_n)^{2g}}{N_1(\log\log\dds_n)^{2g-1}}\delta_n(Q),\]
\noindent en particulier, $\dd_n(Q_1)<\dd_n(Q)$.
\end{cor}
\demo Par d\'efinition de l'indice d'obstruction, on a $\dd_n(Q_1)\leq (\deg_{\L_{K_n}}V)^{\frac{1}{\textnormal{codim}(V)}}$. Ainsi en appliquant la proposition pr\'ec\'edente, on en d\'eduit
\[\dd_n(Q_1)\leq  c_{15}\frac{C_0^{2g+\frac{3}{2}+\rho(r!-1)}(\log\dds_n)^{2g+\rho(r!-1)}}{N_r^{\frac{1}{g-d}}(\log\log\dds_n)^{2g+\rho(r!-1)}}(\log\log\dds_n)^{\frac{1}{g-d}}\delta_n(Q).\]
\noindent \`A partir de l\`a le calcul se fait exactement comme celui de la scolie 7.2. de \cite{davidhindry}\hfill$\Box$

\subsection{Descente finale et preuve du th\'eor\`eme principal \ref{semi}}
\noindent On va maintenant montrer le th\'eor\`eme. Pour cela on choisit $P$ un point de $A(\overline{K})$ et on suppose par l'absurde qu'il ne v\'erifie pas la conclusion du th\'eor\`eme \ref{semi} : on suppose donc que
\[\h_{\L}(P)	 < \frac{c}{\dd_n(P)}\left(\frac{\log\log \dds_n}{C_0\log\dds_n}\right)^{(g+1)!\rho_{\max}},\]
\noindent et on suppose \'egalement que $P$ n'est contenue dans aucune sous-$K_n$-vari\'et\'e de torsion, $B$, stricte de $A_{K_n}$ telle que
\[\left(\deg_{\L_{K_n}} B\right)^{\frac{1}{\text{codim }B}}\leq \delta_n(P)\left(C_0\log\dds_n\right)^{2g+2(g+2)(g+1)!(2g.g!)^g(2g+5)}.\]

\noindent Le point $P$ v\'erifie l'hypoth\`ese $(H)$ pour la valeur de $\rho$ maximale, not\'ee $\rho_{\max}$. La descente que l'on va maintenant expliquer est tout \`a fait similaire au paragraphe 7.3. de \cite{davidhindry}, \`a ceci pr\`et que l'on travaille sur le corps $K_n$ et surtout que l'on ne travaille qu'avec les indices d'obstructions.

\medskip

\noindent On d\'efinit $g$ ensembles de premiers correspondant aux isog\'enies de Frobenius donn\'ees par les ensembles $\P_i^{(j)}$, pour $i\in\{1,\ldots,g\}$ et $j\in\{1,\ldots,g\}$ :
\[\P_i^{(j)}=\{\Id\}\cup\left\{\alpha_v\ |\ \ \frac{N_i^{(j)}}{2}\leq N(v)\leq N_i^{(j)},\ \ v\mid p\text{ non-ramifi\'e dans }K_n\right\},\]
\noindent o\`u les $N_i^{(j)}$ sont d\'efinis par la formule :
\[ N_i^{(j)}=\left(\frac{C_0\log \dds_n}{\log\log\dds_n}\right)^{i.i!\rho_j},\]
\noindent avec $\rho_j=(2gg!)^{g-j}\rho_{\text{min}}$.

\medskip

\noindent Nous  utiliserons aussi des ensembles exceptionnels $\mathcal{S}_i^{(j)}$ dont le cardinal est au plus la moiti\'e du cardinal de $\P_i^{(j)}$.

\medskip

\noindent On introduit deux autres familles d'ensembles d'isog\'enies :
\[\forall j\in \{1,\ldots,g\}\hspace{.7cm}\mathcal{Q}_i=\left\{\b_j\ |\ \ \b_j=\a_g^{(j)}\circ\cdots\circ\a_2^{(j)},\ \ \alpha_i^{(j)}\in\P_i^{(j)}\right\}.\]
\noindent En appliquant la proposition \ref{resume} \`a un point $Q$, pour les ensembles $\P_i^{(j)}$, on obtient un point $Q_1=\b_j(Q)$ avec $\b_j\in \mathcal{Q}_j$. On pose \'egalement $\mathcal{R}_0=\{\textnormal{Id}\}$ et 
\[\forall i\in \{1,\ldots,g\},\hspace{.5cm}\mathcal{R}_i=\left\{F_i\ |\ \ F_i=\b_i\circ\cdots\circ\b_1,\ \ \b_i\in\mathcal{Q}_i\right\}.\]

\medskip

\noindent \textbf{Notations.} De m\^eme que pour les isog\'enies $\alpha$, on associe aux isog\'enies $\b$ et $F$ les op\'erateurs $\bb$ et $\tF$ d\'efinis en rempla\c{c}ant les $\alpha$ par les $\aa$.

\medskip

\noindent On montre maintenant que l'on peut, partant de $P$, extrapoler $g$ fois.

\medskip

\begin{lemme}\label{hypoh}Soient $i\in\{0,\ldots,g-1\}$ un entier et $F_i$ un \'el\'ement de $\mathcal{R}_i$. Le point $P_i=\tF_i(P)$ v\'erifie l'hypoth\`ese (H) avec $\rho=\rho_{i+1}$.
\end{lemme}
\demo Par construction des ensembles $\P_i^{(j)}$ et par choix des $N_i^{(j)}$, l'isog\'enie $F_i$ v\'erifie bien le point 2. de l'hypoth\`ese (H). Ainsi il suffit de v\'erifier l'in\'egalit\'e sur la hauteur de $P_i$ pour pouvoir conclure. Or on a
\[\h_{\L}(P_i)=\textnormal{q}(F_i)\h_{\L}(P).\]
\noindent En posant $q_i=\textnormal{q}(F_i)$, on constate donc qu'il suffit de majorer convenablement $q_i$. Par d\'efinition des ensembles $\mathcal{R}_i$, on a
\begin{align*}
q_i	& \leq \prod_{l=1}^i\max\left\{\textnormal{q}(\b_l)\ |\ \b_l\in\mathcal{Q}_l\right\} \leq \prod_{l=1}^i\prod_{k=2}^g N_k^{(l)}\\
	& \leq \left(\frac{C_0\log \dds_n}{\log\log\dds_n}\right)^{(g+1)!(2g.g!)^g\rho_{\text{min}}}
\end{align*}
\noindent On obtient donc, en appliquant le lemme \ref{deltareduit} (plus pr\'ecis\'ement en appliquant l'in\'egalit\'e (\ref{ipl}) de sa preuve),
\begin{align*}
\h_{\L}(P_i)	& \leq \frac{c}{\dd_n(P)}\left(\frac{\log\log \dds_n}{C_0\log\dds_n}\right)^{(g+1)!\rho_{\max}-(g+1)!(2g.g!)^g\rho_{\text{min}}}\\
		& \leq  \frac{c}{\dd_n(P_i)}\left(\frac{\log \log\dds_n}{C_0\log\dds_n}\right)^{(g+1)!\rho_{\max}-(g+1)!(2g.g!)^g\rho_{\text{min}}}q_i^{g}\\
		& \leq  \frac{c}{\dd_n(P_i)}\left(\frac{\log \log\dds_n}{C_0\log\dds_n}\right)^{(g+1)!\left(\rho_{\max}-(g+1)(2g.g!)^g\rho_{\text{min}}\right)}.
\end{align*}
\noindent On a choisi $\rho_{\text{max}}=(g+2)(2g.g!)^g\rho_{\min}$, ce qui conclut.\hfill$\Box$

\medskip

\noindent On passe maintenant \`a la proposition cruciale, qui va nous permettre d'effectuer la descente. Pour la commodit\'e du lecteur, nous conservons la num\'erotation de \cite{davidhindry}.

\medskip

\begin{prop}\label{final} Il existe un entier $k\in\{1,\ldots,g\}$, une suite $V_0,\ldots,V_k$ de sous-$K_n$-vari\'et\'es strictes de $A_{K_n}$ et une suite d'\'el\'ements $F_i=\b_i\circ\cdots\circ \b_1\in\mathcal{R}_i$, $i\in\{1,\ldots,k\}$ (et $F_0=\Id$) v\'erifiant les propri\'et\'es suivantes :
\begin{enumerate}
\item Les dimensions $d_i$ de $V_i$ sont croissantes.
\item La vari\'et\'e $V_i$ passe par $P_i=\tF_i(P)$.
\item La vari\'et\'e $\bb_i^{-1}V_i$ contient $V_{i-1}$ si $i\in\{1,\ldots,k\}$.
\item On a la majoration
\[\Delta_i:=\deg_{\L_{K_n}}V_i\leq  \frac{c_{15}}{N_1^{(i)}}\left(\frac{\dd_n(P_{i-1})C_0^{2g+\frac{3}{2}+\rho_{i+1}(g!-1)}(\log\dds_n)^{2g+\rho_{i+1}(g!-1)}}{(\log\log\dds_n)^{2g-1+\rho_{i+1}(g!-1)}}\right)^{g-d_i}.\]
\item Si $d_i=d_{i+1}$ et si pour tout $j<i\leq k-1$, $d_j<d_{j+1}$, alors le nombre $\textnormal{q}(\b_{i+1})$ est premier avec $\mid G_{V_i}:G_{V_i}^0\mid$.
\item Il existe $l\in\{0,\ldots,k-1\}$ tel que $d_l=d_{l+1}$.
\end{enumerate}
\end{prop}

\medskip

\noindent Cette proposition se prouve en deux \'etapes, tout comme dans l'article de David et Hindry.

\subsubsection{Proposition \ref{final}, premi\`ere \'etape}

\defi Soient $\mathcal{N}$ un ensemble et $\mathcal{S}$ un sous-ensemble. On dit que $\mathcal{S}$ est \textit{exceptionnel (pour $\mathcal{N}$)} si le cardinal de $\mathcal{S}$ est au plus la moiti\'e du cardinal de $\mathcal{N}$.

\medskip

\begin{lemme}\label{etape1} Soit $l\in\{0,\ldots,g-1\}$. Supposons donn\'ee une suite $(W_i^{(l)})_{0\leq i\leq l}$ de sous-$K_n$-vari\'et\'es strictes de $A_{K_n}$ et une suite d'\'el\'ements $F_i=\b_i\circ\cdots\circ \b_1\in\mathcal{R}_i$, $i\in\{1,\ldots,l\}$ (et $F_0=\Id$) v\'erifiant les propri\'et\'es suivantes :
\begin{enumerate}
\item Les dimensions $d_i^{(l)}$ de $W_i^{(l)}$ sont croissantes en $i$ \`a $l$ fix\'e.
\item La vari\'et\'e $W_i^{(l)}$ passe par $P_i=\tF_i(P)$.
\item La vari\'et\'e $\bb_i^{-1}W_i^{(l)}$ contient $W_{i-1}^{(l)}$ si $i\in\{1,\ldots,l\}$.
\item On a la majoration
\begin{equation}\label{wn}
\Delta_i^{(l)}:=\deg_{\L_{K_n}}W_i^{(l)}\leq  \frac{c_{15}}{N_1^{(i)}}\left(\frac{\dd_n(P_{i-1})C_0^{2g+\frac{3}{2}+\rho_{i+1}(g!-1)}(\log\dds_n)^{2g+\rho_{i+1}(g!-1)}}{(\log\log\dds_n)^{2g-1+\rho_{i+1}(g!-1)}}\right)^{g-d_i^{(l)}}.
\end{equation}
\end{enumerate}
Dans ces conditions, il existe un \'el\'ement $\bb_{l+1}(Q_{l+1})$ tel que :
\begin{enumerate}
\item[{5.}] Le nombre $\textnormal{q}(\b_{l+1})$ est premier avec $\mid G_{W_l^{(l)}}:G_{W_l^{(l)}}^0\mid$ et il existe une suite de sous-$K_n$- vari\'et\'es stricte de $A$, $(W_i^{(l+1)})_{0\leq i\leq l+1}$ v\'erifiant les propri\'et\'es pr\'ec\'edentes 1.,2.,3.,4. avec $l$ remplac\'e par $l+1$ et $F_{l+1}=\b_{l+1}\circ F_l$. De plus cette suite v\'erifie la propri\'et\'e :
\item[{6.}] Pour tout $i\in\{0,\ldots,l\}$, on a $W_i^{(l)}\subset W_i^{(l+1)}$.
\end{enumerate}
\end{lemme}
\demo On va appliquer la proposition \ref{resume} en partant :
\begin{enumerate} 
\item Du point $\tF_l(P)$ ;
\item Des ensembles $\P_1^{(l+1)},\ldots,\P_g^{(l+1)}$ o\`u l'on choisit comme ensembles exceptionnels les
\[\mathcal{S}_i^{(l+1)}=\left\{v\in\P_i^{(l+1)}\ | \ \text{pgcd}\left(N(v), \mid G_{W_l^{(l)}}:G^0_{W_l^{(l)}}\mid\right)\not=1\right\}.\]
\end{enumerate}
\noindent Comme rappel\'e dans le lemme \ref{nombrepremier}, le nombre de premiers divisant un entier positif $n$ est au plus polyn\^omial en $\log n$ et le cardinal de la partie discr\`ete des stabilisateurs de $W_l^{(l)}$ est au plus polynomial en le degr\'e de $W_l^{(l)}$ d'apr\`es le point 3. du lemme \ref{degre}. Ainsi en utilisant la propri\'et\'e 4., on a
\[\text{Card }\mathcal{S}_i^{(j)}\leq C_0\log\dds_n.\]
\noindent Comme tout les $N_i^{(j)}$ sont de cardinal au moins $C_0(\log\dds_n)^2$, on est bien assur\'e que les $\mathcal{S}_i^{(l+1)}$ sont exceptionnels. Par ailleurs, le lemme \ref{hypoh} nous assure que $P_i$ v\'erifie l'hypoth\`ese (H). On va donc pouvoir lui appliquer la proposition \ref{resume} avec $\rho=\rho_i$.

\medskip

\noindent Cette proposition nous fournit un \'el\'ement $\b_{l+1}\in\mathcal{Q}_{l+1}$ tel que la propri\'et\'e 5. soit satisfaites (par le choix m\^eme des ensembles exceptionnels). Par ailleurs, on obtient ainsi une sous-$K_n$-vari\'et\'e, $W_{l+1}^{(l+1)}$, stricte de $A$, $K_n$ irr\'eductible, de dimension $d_{l+1}^{(l+1)}$, passant par $P_{l+1}=\bb_{l+1}(P_l)$. Si $W_{l+1}^{(l+1)}$ \'etait de torsion, alors la vari\'et\'e $\tilde{F}_{l+1}^{-1}W_{l+1}^{(l+1)}$ serait \'egalement de torsion. Or cette derni\`ere contient le point $P$, et est de degr\'e au plus
\begin{align*}
\left(\deg_{\L_{K_n}}F_{l+1}^{-1}W_{l+1}^{(l+1)}\right)^{\frac{1}{\text{codim }W_{l+1}^{(l+1)}}} & \leq \textnormal{q}(F_{l+1})\left(\deg_{\L_{K_n}}W_{l+1}^{(l+1)}\right)^{\frac{1}{\text{codim }W_{l+1}^{(l+1)}}}\\
	& \leq  \textnormal{q}(F_{l+1})\frac{LN^2N_1\times\cdots\times N_{g-1}}{T_g},\\
\end{align*}
\noindent la derni\`ere in\'egalit\'e venant de ce que la vari\'et\'e $W_{l+1}^{(l+1)}$ est donn\'ee par le lemme de z\'ero. En rempla\c{c}ant les param\`etres par leur valeur, on en d\'eduit une contradiction, car $P$ n'est contenu par hypoth\`ese dans aucune sous-$K_n$-vari\'et\'e de torsion stricte de degr\'e au plus
\[\left(\deg_{\L_{K_n}} B\right)^{\frac{1}{\text{codim }B}}\leq \delta_n(P)\left(C_0\log\dds_n\right)^{2g+\frac{1}{2}+2(g+1)!(2g.g!)^g(2g+5)}.\]

\medskip 

\noindent Finalement la vari\'et\'e $W_{l+1}^{(l+1)}$ n'est pas de torsion, et on peut appliquer la proposition \ref{resume} pour majorer plus finement son degr\'e :
\begin{align*}
\deg_{\L_{K_n}}W_{l+1}^{(l+1)}	& \leq  c_{15}\dd_n(P_l)^{g-d_{l+1}^{(l+1)}}\frac{C_0^{(2g+\frac{3}{2}+\rho_{l+1}(r_{l+1}!-1))(g-d_{l+1}^{(l+1)})}(\log\dds_n)^{(2g+\rho(r_{l+1}!-1))(g-d_{l+1}^{(l+1)})}}{N_{r_{l+1}}(\log\log\dds_n)^{(2g+\rho_{l+1}(r_{l+1}!-1))(g-d_{l+1}^{(l+1)})-1}}\\
				& \leq  c_{15}\frac{\dd_n(P_l)^{g-d_{l+1}^{(l+1)}}}{N_1^{(l+1)}}\times \left(\frac{C_0^{2g+\frac{3}{2}}(\log\dds_n)^{2g}}{\log\log\dds_n)^{2g-1}}\right)^{g-d_{l+1}^{(l+1)}}.
\end{align*}
\noindent La derni\`ere in\'egalit\'e s'obtient tout comme l'in\'egalit\'e du corollaire \ref{cor2} en suivant le calcul de la preuve de la scolie 7.2 de \cite{davidhindry}. La vari\'et\'e ainsi construite v\'erifie bien l'in\'egalit\'e \ref{wn}.

\medskip

\noindent Comme dans \cite{davidhindry} nous construisons maintenant la vari\'et\'e $W_0^{(l+1)}$. On va pour cela couper $W_0^{(l)}$ par des formes d\'efinissant incompl\`etement $W_{l+1}^{(l+1)}$, en utilisant l'information sur la multiplicit\'e contenue dans le lemme de z\'eros : on coupe $W_0^{(l)}$ par le nombre minimal de formes $G_1,\ldots,G_u$, choisies parmies les formes d\'efinissant incompl\`etement $W_{l+1}^{(l+1)}$ avec multiplicit\'e sup\'erieure \`a $T_g$, donn\'ees par la proposition \ref{resume}, nulles sur $W_{l+1}^{(l+1)}$, tir\'ees en arri\`ere par $\tF_{l+1}^{-1}$ de sorte que 
\[W_0^{(l)}\cdot \tF_{l+1}^{-1}\mathcal{Z}(G_1)\cdots \tF_{l+1}^{-1}\mathcal{Z}(G_u)\]
\noindent a la m\^eme dimension au point $P$ que
\[W_0^{(l)}\cdot \tF_{l+1}^{-1}(W_{l+1}^{(l+1)}).\]
\noindent De plus la proposition \ref{resume} nous assure qu'il existe de telles formes, nulles sur $W_{l+1}^{(l+1)}$ avec multiplicit\'e sup\'erieure \`a $\frac{1}{g}T_g$, de degr\'e au plus $2L^2N_1^{(l+1)}\times\cdots\times N_{r_{l+1}-1}^{(l+1)}$. Ainsi la vari\'et\'e $W_0^{(l+1)}$ peut \^etre d\'efinie par r\'ecurrence sur $u$. Parmi les composantes isol\'ees de \mbox{$W_0^{(l)}\cdot \tF_{l+1}^{-1}\mathcal{Z}(G_l)$}, on en choisit une, que l'on note $W_{1}'^{(l+1)},$ contenant une composante isol\'ee de dimension maximale de 
\[W_0^{(l)}\cdot \tF_{l+1}^{-1}W_{l+1}^{(l+1)},\]
\noindent passant par le point $P$. On applique maintenant le lemme \ref{bezout} et on obtient
\[\deg_{\L_{K_n}}W_1'^{l+1}\leq \frac{c_{15}}{T_g}\left(\deg_{\L_{K_n}}W_0^{(l)}\right)L^2N_1^{(l+1)}\times\cdots\times N_{r_{l+1}-1}^{(l+1)}.\]
\noindent Par r\'ecurrence sur $u$ on obtient une sous-$K_n$-vari\'et\'e stricte de $A$, $W_u'^{(l+1)}=:W_0^{(l+1)}$, irr\'eductible sur $K_n$, contenant le point $P$ et v\'erifiant
\[W_0^{(l+1)}\subset \tF_{l+1}^{-1}W_{l+1}^{(l+1)},\]
\noindent et de degr\'e major\'e par
\[\deg_{\L_{K_n}} W_0^{(l+1)}\leq \frac{c_{17}}{T_g^u}\deg_{\L_{K_n}}W_0^{(l)}\left(L^2N_1^{(l+1)}\times N_{r_{l+1}-1}^{(l+1)}\right)^u.\]
\noindent En utilisant le fait que $\rho_{l+1}\leq \rho_1$ et $r_{l+1}\leq g$, on remplace maintenant les param\`etres par leur valeur pour obtenir :
\[\deg_{\L_{K_n}}W_0^{(l+1)}\leq c_{17}\dd_n(P_l)^u\left(\frac{C_0^{2g+\frac{3}{2}+\rho_{1}(g!-1)}(\log\dds_n)^{2g+\rho_{1}(g!-1)}}{(\log\log\dds_n)^{2g+\rho_{1}(g!-1)}}\right)^u.\]

\medskip

\noindent En rempla\c{c}ant $\deg_{\L_{K_n}}W_0^{(l)}$ par sa majoration donn\'ee par l'hypoth\`ese de r\'ecurrence et en notant que $u=\dim W_0^{(l)}-\dim W_0^{(l+1)}$, on obtient l'in\'egalit\'e
\[\deg_{\L_{K_n}}W_0^{(l+1)}\leq c_{18}\dd_n(P)^{\text{codim}_A(W_0^{(l+1)})}\left(\frac{C_0^{2g+\frac{3}{2}+\rho_{1}(g!-1)}(\log\dds_n)^{2g+\rho_{1}(g!-1)}}{(\log\log\dds_n)^{2g+\rho_{1}(g!-1)}}\right)^{\text{codim}_A(W_0^{(l+1)})}.\]
\noindent Notons que l'on a ici utilis\'e le corollaire \ref{cor2} pour majorer $\dd_n(P_l)$ par $\dd_n(P)$. Ceci ach\`eve la construction au rang $i=0$. De plus, la propri\'et\'e 6. est bien v\'erifi\'ee pour $i=0$.

\medskip

\noindent On se donne maintenant un entier $m\in\{0,\ldots,l-1\}$ et on suppose les vari\'et\'es $W_i^{(l+1)}$ construites pour $i\in\{0,\ldots,m\}$. On veut construire les vari\'et\'es $W_m^{(l+1)}$ comme pr\'e\-c\'e\-demment. De fait ceci marche effectivement de la m\^eme fa\c{c}on et est d\'etaill\'e dans \cite{davidhindry} p.66-67. \hfill$\Box$

\medskip

\noindent Le lemme \ref{etape1} nous permet d'obtenir le r\'esultat suivant :

\medskip

\begin{lemme}\label{etape11}Soit $u$ un entier compris entre $0$ et $g-1$. S'il existe un \'el\'ement de $F_u=\b_u\circ\cdots\circ\b_1\in\mathcal{R}_u$ et une suite de sous-$K_n$-vari\'et\'es $W_i^{(u)}$ strictes de $A$, v\'erifiant les hypoth\`eses du lemme \ref{etape1}, alors, il existe un \'el\'ement
\[F_g=\b_g\circ\cdots\circ\b_{u+1}\circ F_u\in\mathcal{R}_g\]
\noindent et des sous-$K_n$-vari\'et\'es strictes $(W_i^{(j)})_{0\leq i\leq j,u\leq g}$ stricte de $A$ telles que : pour tout $l\in\{u,\ldots,g\}$, la suite $(W_i^{(l)})_{0\leq i\leq l}$ v\'erifie les propori\'et\'es 1., 2., 3. et 4. du lemme \ref{etape1} pour $F_l=\b_l\circ\cdots\circ\b_1$ et telles que de plus les deux propri\'et\'es suivantes soient v\'erifi\'ees :
\begin{enumerate}
\item[{5'.}] Pour tout $l\in \{u,\ldots,g-1\}$, le nombre $\textnormal{q}(\b_{l-1})$ est premier avec $\mid G_{W_l^{(l)}}:G^0_{W_l^{(l)}}\mid$.
\item[{6'.}] Pour tout $i\in\{0,\ldots,g\}$ et pour tout $l\in\{u,\ldots,g-1\}$, on a
\[W_i^{(l)}\subset W_i^{(l+1)}.\]
\end{enumerate}
\noindent De plus, pour $F_0=\Id$, il existe une telle sous-vari\'et\'e $W_0^{(0)}$.
\end{lemme}
\demo Soient $W_0^{(0)}$ une sous-$K_n$-vari\'et\'e stricte de $A$ passant par $P_0=P$, de dimension minimale et r\'ealisant $\dd_n(P)$, \textit{i.e.}, telle que 
\[\dd_n(P)^{\text{codim }W_0^{(0)}}=\deg_{\L_{K_n}}W_0^{(0)}.\]
\noindent On note $d_0^{(0)}$ sa dimension et on pose $F_0=\Id$. La vari\'et\'e $W_0^{(0)}$ v\'erifie les propri\'et\'es 2. et 4. du lemme \ref{etape1}. De plus, pour $l=0$, les deux autres conditions sont vides, donc v\'erifi\'ees. On applique le lemme \ref{etape1} et on obtient des vari\'et\'es $W_0^{(1)}$, $W_1^{(1)}$ ainsi qu'un \'el\'ement $F_1\in\mathcal{R}_1$. Par r\'ecurrence sur $l$ on obtient alors le lemme. De m\^eme lorsque l'on part d'un entier positif $u$ quelconque, la m\^eme r\'ecurrence permet de conclure.\hfill$\Box$

\subsubsection{Proposition \ref{final}, seconde \'etape}
\noindent Cette seconde \'etape, qui est purement combinatoire, se reprend mot pour mot dans la seconde \'etape de la descente de \cite{davidhindry} p.68-71. Pour ne pas alourdir ce papier plus que de raison, nous l'omettons ici. Cette \'etape permet de prouver la proposition \ref{final}.

\subsection{Conclusion}
\noindent Soit $i$ le plus petit entier compris entre $0$ et $k$ tel que $d_i=d_{i+1}$ dans la proposition \ref{final}. La propri\'et\'e 3. de cette proposition nous assure que $V_i$ est une composante isol\'ee de $\bb_{i+1}^{-1}V_{i+1}$. De plus, la vari\'et\'e $\bb_{i+1}^{-1}V_{i+1}$ est stable par translation par les points de $\text{Ker }\b_{i+1}.$ Ainsi pour tout \'el\'ement $\xi$ de ce noyau, la vari\'et\'e $V_i+\xi$ est une composante isol\'ee de $\bb_{i+1}^{-1}V_{i+1}$. Par la propri\'et\'e 5. de la proposition \ref{final}, le nombre $\textnormal{q}(\b_{i+1})$ est premier avec le cardinal de la partie discr\`ete du stabilisateur de $V_i$. Ainsi, en notant $s_i$ la dimension de ce stabilisateur, on en d\'eduit que le nombre de telles composantes est 
\[\textnormal{q}\left(\b_{i+1}\right)^{g-s_i}.\]
\noindent En comparant les degr\'es, on obtient l'in\'egalit\'e
\[\textnormal{q}(\b_{i+1})\deg_{\L_{K_n}} V_i\leq\deg_{\L_{K_n}}\left(\b_{i+1}^{-1}V_{i+1}\right).\]
\noindent La vari\'et\'e $V_i$ passe par construction par le point $P_i$, donc son degr\'e sur $K_n$ est minor\'e par
\[\deg_{\L_{K_n}}V_i\geq\dd_n(P_i)^{g-d_i}.\]
\noindent La propri\'et\'e 4. de la proposition \ref{final} (majoration de $\Delta_i$) nous donne alors
\[\textnormal{q}(\b_{i+1})^{g-s_i}\delta(P_i)^{g-d_i}\leq \textnormal{q}(\b_{i+1})^{g-d_i}\frac{\dd_n(P_i)^{g-d_i}}{N_1^{(i+1)}}\left(C_0^{2g+\frac{3}{2}+\rho_{i+2}(g!-1)}(\log\dds_n)^{2g+\rho_{i+2}(g!-1)}\right)^{g-d_i}.\]
\noindent En simplifiant par $\dd_n(P_i)$ et en rempla\c{c}ant $N_1^{(i)}$ par sa valeur, on obtient
\[\frac{\left(C_0^{2g+\frac{3}{2}+\rho_{i+2}(g!-1)}(\log\dds_n)^{2g+\rho_{i+2}(g!-1)}\right)^{g-d_i}(\log\log\dds_n)^{\rho_{i+1}}}{C_0^{\rho_{i+1}}(\log\dds_n)^{\rho_{i+1}}}\geq 1.\]

\noindent Ainsi, on en d\'eduit une contradiction si
\[ \left(2g+\frac{3}{2}+\rho_{i+2}(g!-1)\right)(g-d_i)<\rho_{i+1}.\]
\noindent Or, par construction on a
\[\rho_{i+1}\geq (2g.g!)\rho_{i+2}.\]
\noindent Ainsi, si $\rho_{\text{min}}=2g+5$ on peut conclure. Le choix de $\rho_{\text{min}}$ nous permet donc de finir la preuve.\hfill$\Box$

\section{Preuve du th\'eor\`eme \ref{ray}\label{pararay}}

\noindent La preuve repose essentiellement sur trois points : 
\begin{enumerate}
\item En relisant la preuve de R\'emond, on peut dans le cas C.M. utiliser une estimation du cardinal des points de torsion meilleure que celle qu'il utilise : l\`a o\`u il utilise une estimation de Masser, on peut dans notre cas utiliser le corollaire 1.2 de \cite{matorsion}.
\item Dans une vari\'et\'e ab\'elienne $A^n$, avec $A$ simple, les sous-groupes alg\'ebriques sont de dimension un multiple de la dimension de $A$.
\item Notre r\'esultat (th\'eor\`eme \ref{semi}) sur le probl\`eme de Lehmer permet de gagner $1$ dans l'estimation finale.
\end{enumerate}

\subsection{Pr\'eliminaires}

\noindent Notons tout d'abord que le probl\`eme que l'on consid\`ere est stable par isog\'enies. Dans la suite on se restreint donc au cas d'une vari\'et\'e ab\'elienne de type C.M. produit de vari\'et\'es ab\'eliennes g\'eom\'etriquement simples, $A=\prod_{i=1}^m A_i^{n_i}$, les $A_i$ \'etant de dimension $g_i$ et deux \`a deux non isog\`enes. Par ailleurs, on fixe une courbe $X$ qui est transverse dans $A$. R\'emond obtient \'egalement dans le cas C.M., en appliquant la th\'eor\`eme \ref{t1} de David-Hindry, le r\'esultat inconditionnel suivant (c'est le corollaire 1.2., page 529, de \cite{remond})~:

\medskip

\begin{theor}\label{p1}\textnormal{\textbf{(R\'emond)}} L'ensemble $X(\overline{K})\cap A^{[r]}$ est fini d\`es que 
\[r \geq 2+\sum_{i=1}^mg_i.\]
\end{theor}

\medskip

\noindent En relisant la preuve de R\'emond, on constate en fait qu'il prouve un r\'esultat un peu plus fin. On donne pour cela une notation en suivant \cite{matorsion} :

\medskip

\noindent \textbf{Notation.} Soit $A/K_0$ une vari\'et\'e ab\'elienne, on note 
\[\gamma(A)=\inf\left\{b>0\ |\ \exists C(A)>0\ \forall K/K_0 \text{ finie },\left|(A(K)_{\textnormal{tors}})\right|\leq C(A/K_0)[K:K_0]^{b}\right\}.\]

\medskip

\begin{propo}\label{p2} \textnormal{\textbf{(R\'emond)}} L'ensemble $X(\overline{K})\cap A^{[r]}$ est fini d\`es que 
\[r > 1+\sum_{i=1}^m\gamma(A_i).\]
\end{propo}

\medskip

\noindent Dans sa preuve du th\'eor\`eme \ref{p1} ci-dessus, il utilise la majoration due \`a Masser \cite{masser}
\[\textnormal{Card\,}(A_i(K)_{\textnormal{tors}})\ll D^{g_i+\varepsilon}\]
\noindent avec $\varepsilon$ assez petit. En prenant le produit, on voit que $\gamma(A)\leq\sum_{i=1}^mg_i+\varepsilon$ d'o\`u, vue la proposition \ref{p2}, le choix de $r$ dans le th\'eor\`eme \ref{p1}. Pr\'ecis\'ement, dans son article l'utilisation du r\'esultat de Masser est faite dans le corollaire 5.1., page 540, de \cite{remond}.

\medskip

\noindent Il y a donc deux mani\`eres de raffiner ce r\'esultat. La premi\`ere consiste \`a remplacer le terme $1+\sum\gamma(A_i)$ par un terme plus petit. Le th\'eor\`eme principal de notre article (le th\'eor\`eme \ref{semi}) nous permet pr\'ecis\'ement d'am\'eliorer ceci (en rempla\c{c}ant $1+\sum\gamma(A_i)$ par $\sum\gamma(A_i)$). C'est l'objet de la proposition \ref{p11} ci-apr\`es . La seconde am\'elioration possible consiste \`a obtenir une majoration plus pr\'ecise que celle de Masser pour l'exposant $\gamma(A)$ dans le cas d'une vari\'et\'e ab\'elienne simple de type C.M. C'est ensuite la conjonction de ces deux am\'eliorations ainsi qu'une remarque qui nous permettra de prouver le th\'eor\`eme \ref{ray}. Nous expliquons ceci dans le paragraphe \ref{paraconc}. Passons maintenant \`a l'\'enonc\'e et la preuve de la premi\`ere am\'elioration.

\subsection{La premi\`ere am\'elioration}
\noindent On commence par donner l'\'enonc\'e, et on consacre le reste du paragraphe \`a sa preuve. Pour simplifier la v\'erification au lecteur, on s'efforce de conserver les notations de \cite{remond}. Ainsi on notera dans ce qui suit, $r'$ ce que l'on notait $r$ pr\'ec\'edemment. Par ailleurs on note $K_0$ le corps de base, $K=K_0(P)$ une extension de degr\'e $D=[K:K_0]$ de $K_0$ (cf. \cite{remond} page 538--539). 

\medskip

\begin{prop}\label{p11} Soient $A/K_0$ une vari\'et\'e ab\'elienne de type C.M. et $X$ une courbe transverse dans $A$. L'ensemble $X(\overline{K_0})\cap A^{[r']}$ est fini d\`es que 
\[r' > \sum_{i=1}^m\gamma(A_i).\]
\end{prop}

\medskip

\noindent On introduit \'egalement la notation $K_n=K_0(A[n])$ o\`u $n$ est le plus grand ordre des points de torsion de $A(K)$, et on pose $D_{n}=[K_n(P):K_n]$. Notons $T_n$ un point d'ordre $n$ de $A(K)$. On a le diagramme suivant :

$$
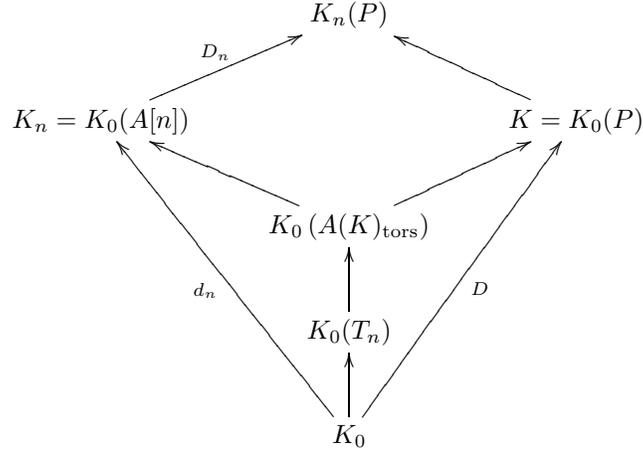
\xymatrix{
 			& K_n(P)					&			\\
K_n=K_0(A[n]) \ar[ur]^{D_n}	&						& K=K_0(P)\ar[ul]	\\
			& K_0\left(A(K)_{\tors}\right)\ar[ur]\ar[ul]	&			\\
			& K_0(T_n)\ar[u]				& 			\\
			& K_0	 \ar[u] \ar[uuur]_D\ar[uuul]^{d_n}	&  	}
$$

\noindent Rappelons un lemme classique dont nous avons besoin.

\medskip

\begin{lemme}\label{cori} Soient $n$ un entier positif et $A/K_0$ une vari\'et\'e ab\'elienne de dimension $g$. On a
\[ d_n:=[K_0(A[n]):K_0] \leq n^{4g^2}.\]
\end{lemme}
\demo Soit $n\geq 1$ un entier. La repr\'esentation naturelle 
\[\rho : \Gal(\overline{K_0}/K_0)\rightarrow \text{Aut}\left(A[n]\right)\]
\noindent nous donne une injection de $\Gal(K_0(A[n])/K_0)$ dans $\text{GL}_{2g}(\Z/n\Z)$. Ceci conclut. \hfill$\Box$

\medskip

\begin{cor}\label{faible}Il existe deux constantes $c_1$ et $c_2$ strictement positives ne d\'ependant que de $A/K_0$, telles que l'on a l'in\'e\-ga\-li\-t\'e
\[d_n\leq c_1 D^{c_2}.\]
\end{cor}
\demo En utilisant l'in\'egalit\'e du lemme \ref{cori}, on majore $d_n$ par une puissance de $n$. Par ailleurs, on sait (en utilisant par exemple les r\'esultats transcendants de Masser \cite{masser}, ou alg\'ebrique de Silverberg \cite{silverberg} dans notre cas) que l'on peut majorer $n$ par une puissance du degr\'e de $[K_0(T_n):K_0]\leq D$. Ceci permet de conclure.\hfill$\Box$

\medskip

\noindent \textbf{Preuve de la proposition \ref{p11} :} ceci se fait selon les deux \'etapes suivantes~:

\medskip

\noindent \textbf{\'Etape 1 : } On montre que si $r'>\sum_{i=1}^m\gamma(A_i)$ alors on peut majorer $D$ en fonction de $D_n$.

\noindent \textbf{\'Etape 2 : } En reprenant le paragraphe 7, pages 545--547, de \cite{remond} en travaillant sur $K_n$ plut\^ot que sur $K^t$, on montre que $D_n$ est born\'e d\`es que $r'\geq 2$.

\medskip

\noindent Ainsi la conjonction des deux \'etapes entra\^ine que $D$ est born\'e  pour $r'>\sum_{i=1}^m\gamma(A_i)\geq 1$. De plus le lemme 3.3, page 535, de \cite{remond} nous dit que l'ensemble $X(\overline{K})\cap A^{[r']}$ (et m\^eme l'ensemble $X(\overline{K})\cap A^{[1]}$) est de hauteur born\'ee. Le th\'eor\`eme de Northcott nous permet alors de conclure concernant la finitude de $X(\overline{K})\cap A^{[r']} $. Ceci conclut donc la preuve de la proposition \ref{p11} modulo les \'etapes 1 et 2 pr\'ec\'edentes. Notons \`a titre de remarque que nous n'avons pas besoin ici d'appliquer le lemme 7.1 de \cite{remond} utilisant le th\'eor\`eme de Raynaud (ex-conjecture de Manin-Mumford).

\medskip

\noindent Il nous reste, pour compl\'eter la preuve de la proposition \ref{p11}, \`a prouver les deux \'etapes pr\'ec\'edentes. C'est ce qu'on fait dans les deux sous-paragraphes suivants. Dans ces deux \'etapes, on appliquera notre th\'eor\`eme \ref{semi} en direction du probl\`eme de Lehmer relatif, en conjonction avec le corollaire \ref{faible} afin de majorer $d_nD_n$ par une puissance de $D$.

\subsubsection{Rappels de notations de \cite{remond}} 
\noindent On travaille sur une vari\'et\'e ab\'elienne de type C.M., $A=\prod_{i=1}^mA_i^{n_i}$ d\'efinie sur un corps de nombres $K_0$.
On note $X$ la courbe transverse incluse dans $A$. Comme dans \cite{remond} paragraphe 5, p.538, on note  
\[P=(P_{1,1},\ldots,P_{1,n_1},\ldots,P_{i,1},\ldots,P_{i,n_i},\ldots,P_{m,n_m})\]
\noindent le point de degr\'e $D=[K_0(P):K_0]$ sur $K_0$ avec lequel on travaille (dans \cite{remond} il s'agit du point $P'=f(P)$ o\`u $f$ est une isog\'enie fixe entre la vari\'et\'e ab\'elienne ambiante et le produit $\prod A_i^{n_i}$ : dans notre situation on a $f=\text{Id}$). Pour tout $i\in \{1,\ldots,m\}$, nous introduisons de plus deux $\text{End}(A_i)$-modules :

\medskip

\begin{enumerate}
\item $N_i$ est le sous-$\text{End}(A_i)$-module de $\text{Hom}(A,A_i)$ des morphismes nuls en $P$.
\item $\Gamma_i$ est le sous-$\text{End}(A_i)$-module de $A_i(\overline{K})$ engendr\'e par les points $P_{i,1},\ldots,P_{i,n_i}$.
\end{enumerate}

\medskip

\noindent Pour tout $i$, l'espace vectoriel r\'eel $\Gamma_i\otimes\mathbb{R}$ est naturellement muni d'une structure euclidienne en utilisant la hauteur de N\'eron-Tate. On note (suivant \cite{remond} p. 539)
\[\nu_i=\text{Vol}\left(\Gamma_i\otimes  \mathbb{R}/\left(\Gamma_i/(\Gamma_i)_{\tors}\right)\right)\]
\noindent le volume pour cette norme.

\medskip

\noindent Avec les notations de \cite{remond}, on prend $r=2$. Pour tout $i\in\{1,\ldots,m\}$, on note  $s_i$ le rang du $\text{End}(A_i)$-module $\Gamma_i$, et on introduit les deux nombres
\[t=\sum_{i=1}^m (n_i-s_i),\ \ \text{ et }\ \  q=\sum_{i=1}^m g_i(n_i-s_i).\]
\noindent L'entier positif $q$ est introduit dans \cite{remond} proposition 6.1., p.543, et la ligne la pr\'ec\'edant. L'entier positif $t$ est quant \`a lui introduit \`a la deuxi\`eme ligne de la preuve de la proposition 6.1 de \cite{remond}.

\medskip

\noindent Enfin les morphismes $\psi_j$ intervenant dans la suite sont des \'el\'ements de $N_i$ apparaissant dans la preuve de la proposition 6.1. de \cite{remond}.

\subsubsection{Preuve de l'\'etape 1}\label{suivant}

\noindent En utilisant les notations et la preuve de la proposition 6.1 page 543--544 de \cite{remond}, on a 
\begin{equation}\label{key}
D\ll \left(\prod_{j=1}^t\| \psi_j\|^{2\text{rg}\psi_j}\right)^{\frac{1}{q}}.
\end{equation}
(Notons que cette in\'egalit\'e a \'et\'e introduite pour la premi\`ere fois dans ce contexte, dans le cadre de $\mathbb{G}_m^n$, dans l'article \cite{BMZ}).

\medskip 

\noindent D\'etaillons un peu l'obtention de cette in\'egalit\'e (\ref{key})~: Dans \cite{remond} p. 544, R\'emond montre d'abord que $P$ est un point isol\'e de $X\cap V$ ( o\`u $V$ est une vari\'et\'e auxiliaire d\'efinie sur $K_0$, introduite p. 544 ligne 1). Ceci implique que 
\[D\leq \deg X\cap V.\]
\noindent Par le th\'eor\`eme de B\'ezout, et $X$ \'etant fixe, on obtient donc
\[D\leq \deg X\deg V\ll \deg V.\]
\noindent Il reste \`a majorer le degr\'e de $V$ : c'est ce qui est fait dans la seconde partie de la page 544 de \cite{remond}. Pr\'ecis\'ement R\'emond obtient
\[\deg V\ll \left(\prod_{j=1}^t\| \psi_j\|^{2\text{rg}\psi_j}\right)^{\frac{r-1}{q}}.\]
\noindent Or nous avons d\'ej\`a rappel\'e que dans notre situation nous prenons $r=2$ (attention  \`a ne pas confondre $r$ et $r'$ dans les notations de \cite{remond}). Ceci conclut donc la preuve de l'in\'egalit\'e (\ref{key}).

\medskip

\noindent Comme expliqu\'e au d\'ebut de la preuve de la proposition 6.1 de \cite{remond}, la famille de morphismes not\'ee $\varphi_{i,j}$ dans la proposition 5.3 de \cite{remond} n'est autre que la famille $\psi_i$ pour $1\leq i\leq t$. On peut donc maintenant relire la proposition 5.3 de \cite{remond} et on obtient : 
\begin{equation}\label{key2}
\prod_{j=1}^t\| \psi_j\|^{2\text{rg}\psi_j}\ll \prod_{i=1}^m\text{vol}\left(N_i\otimes\R/N_i\right)\ll \prod_{i=1}^m\left(\mid A_i(K)_{\tors}\mid v_i^{-1}\right).
\end{equation}
\noindent Explicitons ceci : le terme de gauche de l'in\'egalit\'e est exactement le terme de gauche de l'in\'egalit\'e de la proposition 5.3 de \cite{remond}. En effet, les $\phi_{i,j}$ sont des morphismes non nuls \`a valeurs dans les vari\'et\'es ab\'eliennes simples $A_i$ de dimension $g_i$. Donc $g_i=\text{rg}\phi_{i,j}$ ce qui donne bien le membre de gauche de l'in\'egalit\'e. Le membre de droite de (\ref{key2}) s'obtient en suivant la preuve de la proposition 5.3 de \cite{remond} : plus pr\'ecis\'ement, dans \cite{remond}, le volume $\text{vol}\left(N_i\otimes\R/N_i\right)$ appara\^it avec un exposant $2g_i/d_i$. Mais nous sommes dans le cas particulier de vari\'et\'es ab\'eliennes de type CM, donc $2g_i/d_i=1$. La derni\`ere in\'egalit\'e correspond \`a la fin de la preuve de la proposition 5.3 de \cite{remond} : on applique ses formules (1) et (2) p. 539 et on utilise que le point $P$ est de hauteur born\'ee.

\medskip

\noindent En utilisant la d\'efinition de $\gamma(A_i)$ on a 
\[\mid A_i(K)_{\tors}\mid\ll D^{\gamma(A_i)+\epsilon},\]
\noindent ceci \'etant valable pour tout $\epsilon$, le signe $\ll$ d\'epenant de $\epsilon$ (nous choisirons \`a la fin $\epsilon$ suffisamment petit). Avec l'in\'egalit\'e (\ref{key2}), ceci nous donne 
\begin{equation}\label{eq1}
D\ll \left[\prod_{i=1}^m\left(\mid A_i(K)_{\tors}\mid v_i^{-1}\right)\right]^{\frac{1}{q}}\ll \left(D^{\sum_{i=1}^m \gamma(A_i)+\e}\prod_{i=1}^m v_i^{-1}\right)^{\frac{1}{q}}.
\end{equation}

\noindent \textbf{Conclusion de la preuve de l'\'etape 1 : application de notre th\'eor\`eme \ref{semi}}.

\medskip

\noindent En appliquant notre th\'eor\`eme \ref{semi}, on peut maintenant relire la proposition 5.2, page 541, de \cite{remond} : grace au th\'eor\`eme \ref{semi} on am\'eliore l'estimation faisant intervenir la minoration de la hauteur des points $Q_{i,j}$. On obtient ainsi (rappelons \`a nouveau que dans notre cas de type CM, $g_i/d_i=1/2$)~:
\begin{equation}\label{key3}
\prod_{i=1}^mv_i^{\frac{1}{2}}\geq c_8 D_n^{-\frac{1}{2}}D^{-\frac{\e}{2}}.
\end{equation}
\noindent D\'etaillons l'obtention de cette in\'egalit\'e : on suit la preuve de la proposition 5.3 de \cite{remond}. Comme il l'explique, son appendice fournit une famille $Q_{i,j}$ telle que 
\[\prod_{i=1}^m\prod_{j=1}^{s_i} \widehat{h}(Q_{i,j})^{\frac{g_i}{2}}\ll \prod_{i=1}^mv_i^{\frac{g_i}{d_i}}=\prod_{i=1}^mv_i^{\frac{1}{2}}.\]
\noindent Ensuite, l\`a o\`u il applique la conjecture de Lehmer relatif, nous appliquons notre th\'eor\`eme \ref{semi}, obtenant ainsi
\[\prod_{i=1}^m\prod_{j=1}^{s_i} \widehat{h}(Q_{i,j})^{\frac{g_i}{2}}\gg D_n^{-\frac{1}{2}}(\log d_nD_n)^{-\kappa(g)}.\]
\noindent Il ne reste plus qu'\`a estimer le terme logarithmique, ce qui se fait pr\'ecis\'ement en appliquant notre corollaire \ref{faible}.

\medskip

\noindent En regroupant les in\'egalit\'es (\ref{key2}) et (\ref{key3}), on obtient finalement que $D$ est major\'e polynomialement en fonction de $D_{n}$ si et seulement si 
\[q >\sum_{i=1}^m \gamma(A_i).\]
\noindent Or par construction de $q$, on sait que $q\geq r'$, donc en prenant $r'> \sum_{i=1}^m \gamma(A_i)$ on a bien la conclusion voulue. ceci ach\`eve la preuve de l'\'etape 1 intervenant dans la d\'emonstration de la proposition \ref{p11}. \hfill$\Box$

\subsubsection{Preuve de l'\'etape 2}

\noindent On reprend ce qui est fait au paragraphe 7 de l'article \cite{remond} de R\'emond, en rempla\c{c}ant $K_{\tors}$ et $D_{\tors}$ par $K_n$ et $D_n$. Tout marche pareil : les points de torsion intervenant sont des points de $A(K)_{\tors}$, donc en particulier d\'efinis sur $K_n$. Ceci entraine que $D_n$ est born\'e d\`es que $r'\geq 2$ et conclut donc la preuve de notre proposition \ref{p11}.

\subsection{Conclusion\label{paraconc}}

\noindent Nous expliquons maintenant comment conclure la preuve du th\'eor\`eme \ref{ray}. Comme indiqu\'e auparavant, nous avons besoin pour conclure d'un raffinement des bornes sur la torsion dans les vari\'et\'es ab\'eliennes de type C.M. Dans cette direction, nous obtenons dans le corollaire 1.2 de \cite{matorsion} le r\'esultat suivant~:

\medskip

\begin{theo}\label{monthtors}\textnormal{\textbf{\cite{matorsion}}} Soit $A/K_0$ une vari\'et\'e ab\'elienne de dimension $g$, simple et de type C.M. On a
\[\gamma(A)\leq \frac{2g}{2+\log_2 g}\]
\noindent o\`u $\log_2$ d\'enote le logarithme en base $2$.
\end{theo}

\medskip

\noindent Si $A$ est (isog\`ene \`a) une puissance d'une courbe elliptique C.M., alors le th\'eor\`eme pr\'ec\'edent et la proposition \ref{p11} permettent de conclure. On suppose d\'esormais que $A$ est (isog\`ene \`a) une puissance d'une vari\'et\'e ab\'elienne simple de type C.M. $A_1$ de dimension $g_1$ sup\'erieure \`a $2$. 

\medskip

\noindent On constate que d\`es que $g_1$ est strictement sup\'erieur \`a $1$, le th\'eor\`eme \ref{monthtors} entra\^ine en particulier $\gamma(A_1)<g_1$ ce qui est meilleur que la borne de Masser. N\'eanmoins en appliquant la proposition \ref{p11}, ceci ne permet \textit{a priori} que d'obtenir la finitude de l'ensemble
\[A^{[g_1]}\cap X(\overline{K_0}).\]

\noindent On utilise donc pour conclure le lemme trivial suivant :

\medskip

\begin{lemme}Soit $G$ un sous-sch\'ema en groupes de $A=\prod_{i=1}^mA_i^{n_i}$, les $A_i$ \'etant de dimensions respectives $g_i$. Si $G$ est strictement inclus dans $A$, alors sa codimension v\'erifie
\[\text{codim}(G)\geq \min_{1\leq i\leq m}g_i.\]
\end{lemme}
\demo Soit $G^0$ la composante connexe de l'identit\'e de $G$. C'est une sous-vari\'et\'e ab\'elienne de $A$. Elle est donc isog\`ene \`a un produit de $A_i^{s_i}$. Comme $G$ est strictement inclus dans $A$, il existe $i\in\{1,\ldots,m\}$ tel que $s_i\leq n_i-1$. Ceci conclut.\hfill$\Box$

\medskip

\begin{cor}\label{idiot}Soit $A=\prod_{i=1}^mA_i^{n_i}$, les $A_i$ \'etant de dimensions respectives $g_i$. On suppose que $A$ est de dimension sup\'erieure \`a $2$. On note $A^{[r]}=\bigcup_{\text{codim}(G)\geq r}G(\overline{K})$, et on note $g_{\min}=\min_{1\leq i\leq m} g_i$. On a 
\[ A^{[2]}=A^{[\max\{g_{\min}, 2\}]}.\]
\end{cor}
\demo C'est \'evident par le lemme pr\'ec\'edent.\hfill$\Box$

\medskip

\noindent On termine maintenant la preuve du th\'eor\`eme \ref{ray}~: soit $A_1/K_0$ une vari\'et\'e ab\'elienne de type C.M. de dimension $g_1\geq 2$. Soit $n\geq 1$ un entier et soit $X$ une courbe transverse dans $A=A_1^n$. On sait par ce qui pr\'ec\`ede que $A^{[g_1]}\cap X(\overline{K_0})$ est fini. Le corollaire \ref{idiot} permet donc de conclure : l'ensemble $A^{[2]}\cap X(\overline{K_0})$ est fini.\hfill$\Box$

\medskip

\rem Notons que dans le th\'eor\`eme \ref{ray}, le cas le plus difficile est le cas o\`u $A_1$ est de dimension $2$. En effet, si $A_1$ est de dimension sup\'erieure \`a $3$, notre th\'eor\`eme \ref{monthtors} permet de conclure sans avoir \`a utiliser le raffinement sur le probl\`eme de Lehmer. Par contre en dimension $2$, on peut voir que $\gamma(A_1)=\frac{4}{3}$ (cf. la proposition 1.1 de \cite{matorsion}) et dans ce cas, l'utilisation de notre th\'eor\`eme \ref{semi} est indispensable.

\section{Preuve du th\'eor\`eme \ref{absolu}}

\subsection{Une petite r\'eduction g\'eom\'etrique}

\begin{lemme}\label{hodge}Soient $X/K$ est une vari\'et\'e projective de dimension $g$ et $\L_1,\ldots,\L_g$ des fibr\'es en droites amples. On a l'in\'egalit\'e
\[(\L_1^g)\ldots(\L_g^g)\leq \left(\L_1\cdot\ldots\cdot\L_g\right)^g.\]
\end{lemme}
\demo Il s'agit d'une g\'en\'eralisation en dimension $g$ d'un r\'esultat bien connu pour les surfaces, d\'ecoulant du th\'eor\`eme de l'indice de Hodge. Cette g\'en\'eralisation est elle m\^eme bien connue des sp\'ecialistes (cf. par exemple l'exercice 6 p.50 de \cite{debarre}). Ceci se prouve par r\'ecurrence sur $g$ en se ramenant en dimension inf\'erieure (jusqu'\`a la dimension $2$) gr\^ace au th\'eor\`eme de Bertini.\hfill$\Box$

\medskip

\noindent On note $N$ la forme quadratique d\'efinie positive sur $\textnormal{End}(A)\otimes\R$ d\'eduite de l'involution de Rosati (correspondant au fibr\'e ample $\L$). On note $\|\cdot\|=\sqrt{N(\cdot)}$ la norme qui s'en d\'eduit.

\medskip

\begin{cor}\label{geom}Soit $\alpha$ une isog\'enie de $A/K$. On note $n$ son degr\'e, et $[n]$ l'isog\'enie correspondante. On a
\[\|[n]\|\leq (\L^g)^{\frac{g-1}{2}}\|\alpha\|^{g}.\]
\end{cor}
\demo Rappelons que par d\'efinition (cf. \cite{mumford} p.192 avec un facteur de renormalisation $2g$), on a
\[N(\alpha)=\frac{1}{\L^g}\left(\alpha^*\L\cdot\L^{g-1}\right),\ \text{ et }\ n=\deg \alpha=\frac{1}{\L^g}(\alpha^*\L)^g.\]
\noindent On applique le lemme \ref{hodge} pr\'ec\'edent avec $\L_1=\alpha^*\L$ et $\L_i=\L$ pour $i\geq 2$. On en d\'eduit
\[n\leq \frac{1}{\L^g}(\alpha^*\L)^g(\L^g)^{g-1}\leq \frac{1}{\L^g}\left(\alpha^*\L\cdot\L^{g-1}\right)^g=(\L^g)^{g-1}N(\alpha)^{g}.\]
\noindent Le degr\'e de l'isog\'enie $[n]$ se calcule explicitement : c'est le poids $\textnormal{q}([n])=n^2$ de cette isog\'enie admissible. Ceci permet de conclure.\hfill$\Box$

\subsection{Preuve du th\'eor\`eme \ref{absolu}}

\noindent Soit $A/K$ une vari\'et\'e ab\'elienne simple de type C.M. de dimension $g$. On note que pour tout entier $n\geq 1$ et pour tout entier positif $r$, on a $K(A^n[r])=K(A[r])$. On donne tout d'abord un r\'esultat de comparaison de degr\'e dont nous aurons besoin dans la suite.

\medskip

\begin{lemme}\label{silv} Il existe une constante strictement positive $c_1$, ne d\'ependant que de $A/K$ et $n$, telle pour tout entiers positifs $m$ et $r$, en notant $K_m=K(A[m])$, $K_{mr}=K(A[mr])$ et $T$ un point de torsion de $A^n$ d'ordre exactement $r$, on a
\[ \deg(K_m(T)/K)\leq \deg(K_{mr}/K) \leq c_1\deg(K_m(T)/K)^{16g^2}.\]
\end{lemme}
\demo Il suffit bien sur de montrer l'in\'egalit\'e de droite. On note $g$ la dimension de $A$ et on distingue pour cela deux cas. Si $r\geq m$, alors
\begin{align*}
\deg(K_m(T)/K)	& \geq \deg(K(T)/K)\geq c_1 r^{\frac{1}{2}}\text{ \hspace{.5cm} d'apr\`es le (1.1) de Silverberg \cite{silverberg}}\\
		& \geq c_1(rm)^{\frac{1}{4}} \geq c_1 \deg(K_{mr}/K)^{\frac{1}{16g^2}} \text{ \hspace{1cm} d'apr\`es le lemme \ref{cori}.}
\end{align*}
\noindent Par ailleurs, dans l'autre cas, si $m\geq r$, on a 
\[\deg(K_m(T)/K)\geq\deg(K_m/K)\geq c_1 m^{\frac{1}{2}}\geq c_1 (mr)^{\frac{1}{4}}\geq c_1 \deg(K_{mr}/K)^{\frac{1}{16g^2}}\]
\noindent par le m\^eme argument. Ceci conclut.\hfill$\Box$

\medskip

\noindent On passe maintenant \`a la preuve du th\'eor\`eme \ref{absolu}. Pour cela on raisonne par r\'ecurrence sur la dimension $n$. Si $n=1$, il s'agit du th\'eor\`eme \ref{semi}. On suppose donc le r\'esultat vrai au rang $n-1\geq 1$ et on veut le montrer au rang $n$. On suppose par l'absurde que le r\'esultat est faux en dimension $n$. Ainsi, il existe un point $P=(x_1,\ldots,x_n)\in A^n(\overline{K})$ et une sous-vari\'et\'e $V$ de $A^n$ d\'efinie sur $K_m=K(A[m])$ tels que
\begin{equation}\label{ddd}
\widehat{h}_{\L_n}(P)<C_0(n)^{-1}\ddd^{-1}\left(\frac{\log\log [K_m:K]\ddd}{\log [K_m:K]\ddd}\right)^{\kappa(n)}.
\end{equation}
\noindent Comme le point $P$ est un point $\overline{K}$-rationnel de $V$, on a n\'ecessairement $\dd_{\L,K_m}(P)\leq\ddd$. Ainsi l'in\'egalit\'e \ref{ddd} pr\'ec\'edente nous donne :
\begin{equation}\label{ddd2}
\widehat{h}_{\L_n}(P)<C_0(n)^{-1}\dd_{\L,K_m}(P)^{-1}\left(\frac{\log\log [K_m:K]\dd_{\L,K_m}(P)}{\log [K_m:K]\dd_{\L,K_m}(P)}\right)^{\kappa(n)}.
\end{equation}
\noindent On voit ainsi que l'on est bien dans la seconde partie de l'alternative du th\'eor\`eme \ref{semi}. Ainsi, il existe une sous-vari\'et\'e de torsion stricte $B/K_m$ de $A^n$ dont le degr\'e est major\'e par 
\begin{equation}\label{degb}
\deg_{\L}B^{\frac{1}{\text{codim}B}}\leq c(A/K,\L,n)\ddd\left(\log [K_m:K]\ddd\right)^{-2n-2\kappa(n)}.
\end{equation}
\noindent On \'ecrit
\[B_{\overline{K}}=\bigcup_{\sigma\in \Gal(K_m(T)/K_m)}\left(H+\sigma(T)\right),\]
\noindent o\`u $H$ est la composante connexe de l'origine de $B$ et $T$ est un point de torsion de $A^n$ d'ordre un certain entier positif $r$.

\medskip

\noindent On note 
\[\Lambda=\left\{\phi\in\text{Hom}(A^n,A)\ |\ H\subset \text{Ker}\phi\right\}.\]
\noindent Rappelons que l'on note $N$ la forme quadratique d\'efinie positive sur $\textnormal{End}(A)\otimes\R$ d\'eduite de l'involution de Rosati (correspondant au fibr\'e ample $\L$), et $\|\cdot\|=\sqrt{N(\cdot)}$ la norme qui s'en d\'eduit. En utilisant l'isomorphisme $\text{Hom}(A^n,A)\simeq \text{End}(A)^n$, on munit $\text{Hom}(A^n,A)\otimes \mathbb{R}$ de la norme 
\[(f_1,\ldots,f_n)\mapsto \max_{1\leq i\leq n}\|f_i\|\]
\noindent On note encore $\|\cdot\|$ cette norme. On note \'egalement $\text{Vol}(\Lambda)$ le volume, pour la norme pr\'ec\'edente, de $(\Lambda\otimes \R)/\Lambda$. Par le th\'eor\`eme $\hat{2}$ de \cite{bertrand} on sait qu'il existe un \'el\'ement $\phi$ non nul de $\Lambda$ tel que 
\[\|\phi\|\leq c_3(\deg_{\L}H)^{c_4}.\]
\noindent En utilisant le corollaire \ref{geom}, on voit que quitte \`a modifier $c_4$, on peut en fait supposer que $\phi=(m_1,\ldots,m_n)$ o\`u les $m_i$ sont les isog\'enies admissibles ``multiplication par $m_i$''. (Il suffit de remplacer chacune des composante $\phi_i$ de $\phi$ par $m_i=\hat{\phi_i}\circ\phi_i$ o\`u $\hat{\phi_i}$ est l'isog\'enie duale de $\phi_i$). Enfin, quitte \`a renum\'eroter, on peut supposer que 
\[\|\phi\|=m_n.\]
\noindent On pose maintenant $H'$ la composante connexe de l'origine de $\text{ker}\phi$. Par construction $H'$ contient $H$ (car $H$ est une vari\'et\'e ab\'elienne contenue dans $\text{ker}\phi$). De plus $A$ \'etant simple, la vari\'et\'e ab\'elienne $H'$ est param\'etr\'ee par
\[\Phi : A^{n-1}\rightarrow A^n,\hspace{1cm} (x_1,\ldots,x_{n-1})\mapsto \left(m_n(x_1),\ldots,m_n(x_{n-1}),-\sum_{i=1}^{n-1}m_i(x_i)\right).\]
\noindent Quitte \`a remplacer $T$ par $\sigma(T)$, on peut supposer que $P\in H+T$. On se donne maintenant $y\in\Phi^{-1}(P-T)$ et $V'=\Phi^{-1}(V-T)$.

\medskip

\begin{lemme}Avec les notations pr\'ec\'edentes, les propri\'et\'es suivantes sont v\'erifi\'ees :
\begin{enumerate}
\item $V'\subsetneq A^{n-1}$.
\item $V'$ est incompl\`etement d\'efinie sur $K_{mr}$ dans $A^{n-1}$ par des \'equations de degr\'e inf\'erieur \`a $2^n\|\phi\|^2\ddd$.
\item $y\in (V')^{\star}$.
\item $\widehat{h}_{\L_{n}}(P)\geq \|\phi\|^2\widehat{h}_{\L_{n-1}}(y).$
\end{enumerate}
\end{lemme}
\demo Les points 1. et 3. sont faciles. Pour le point 2., les composantes de $\Phi$ \'etant admissibles, on voit (cf. \cite{hindry} lemme 6 (iii)) que $V'$ est incompl\`etement d\'efinie par des \'equations de degr\'e major\'e par
\[ \max\left\{m_n^2,\left\|\sum_{i=1}^{n-1}m_i\right\|^2\right\}\ddd\leq 2^nm_n^2\ddd=2^n\|\phi\|^2\ddd.\]
\noindent La minoration de hauteur d\'ecoule des in\'egalit\'es suivantes :
\begin{align*}
\h_{\L_n}(P)	& =\h_{\L_n}(P-T)=\h_{\L_n}(\Phi(y))=\sum_{i=1}^{n}\h_{\L}(\Phi_i(y))\\
		& \geq \sum_{i=1}^{n-1}\h_{\L}(\Phi_i(y))=\sum_{i=1}^{n-1}\h_{\L}(m_n(y_i))\\
		& \geq m_n^2 \h_{\L_{n-1}}(y)=\|\phi\|^2\h_{\L_{n-1}}(y).
\end{align*}
\noindent Ceci conclut.\hfill$\Box$

\medskip

\noindent Ce lemme nous permet de terminer la preuve par r\'ecurrence : par hypoth\`ese de r\'ecurrence, il existe une constante $c(A/K,\L,n-1)$ strictement positive telle que
\[\|\phi\|^{-2}\h_{\L_{n}}(P)\geq \h_{\L_{n-1}}(y)\geq \frac{c(A/K,\L,n-1)^{-1}}{2^n\|\phi\|^2\ddd}\left(\frac{\log\log\left([K_{mr}:K]n\|\phi\|^2\ddd\right)}{\log\left([K_{mr}:K]n\|\phi\|^2\ddd\right)}\right)^{\kappa(n-1)}.\]
\noindent En utilisant la majoration de $\|\phi\|$ en fonction de $\deg_{\L}H$, en utilisant le lemme \ref{silv} et en utilisant \'egalement l'identit\'e
\[\deg_{\L}B=[K_m(T):K_m]\deg_{\L}H,\]
\noindent on obtient :
\[\h_{\L_{n}}(P)\geq \frac{c(A/K,\L,n-1)^{-1}c_5}{n\ddd}\left(\log\left([K_{m}:K]n\deg_{\L}B\ddd\right)\right)^{-\kappa(n-1)}.\]
\noindent On utilise maintenant l'in\'egalit\'e (\ref{degb}) et le fait que $\kappa(n)-\kappa(n-1)>0$ pour conclure.

\newpage

\appendix
\section{Appendice}

\noindent On montre ici que la premi\`ere partie de la conjecture de Lehmer ab\'elienne \ref{conj1app2} (minoration des points engendrant la vari\'et\'e ab\'elienne en terme de l'indice d'obstruction), formul\'ee dans \cite{davidhindry} entra\^ine la seconde partie de cette conjecture (minoration des points non de torsion en fonction du degr\'e du point et de la dimension du plus petit sous-groupe alg\'ebrique contenant le point). De m\^eme pour le r\'esultat non-conjectural, ce qui permet d'am\'eliorer le pr\'ec\'edent meilleur r\'esultat connu, d\^u \`a Masser \cite{lettre}, pour la minoration des points d'ordre infini sur les vari\'et\'es ab\'eliennes de type C.M. Par ailleurs on montre que la conjecture de Lehmer ab\'elienne entra\^ine la conjecture de Lehmer ab\'elienne multihomog\`ene \textit{a priori} plus forte, telles qu'elles sont \'enonc\'ees dans \cite{davidhindry}. On montre \'egalement que toute avanc\'ee en direction de la conjecture de Lehmer entra\^ine une avanc\'ee similaire en direction de la conjecture multihomog\`ene. En utilisant le r\'esultat principal de \cite{davidhindry} on en d\'eduit, en direction de la conjecture multihomog\`ene, une minoration optimale aux puissances de log pr\`es dans le cas des vari\'et\'es ab\'eliennes de type C.M.

\subsection{Sur la conjecture de Lehmer sur les vari\'et\'es ab\'e\-lien\-nes}

\noindent En utilisant le th\'eor\`eme  de David et Hindry \cite{davidhindry}, on obtient un r\'esultat, optimal aux puissances de $\log$ pr\`es en direction de l'in\'egalit\'e (\ref{e2}) du probl\`eme de Lehmer ab\'elien (conjecture \ref{conj1app2}). 

\medskip

\begin{theo}\label{thamel}Si $A/K$ est de type C.M., alors il existe une constante strictement positive $c(A/K,\mathcal{L})$ telle que pour tout point $P\in A(\overline{K})$ d'ordre infini, on a 
\[\widehat{h}_{\mathcal{L}}(P)\geq\frac{c(A/K,\mathcal{L})}{D^{\frac{1}{g_0}}}\left(\log 2D\right)^{-\kappa(g_0)},\]
\noindent o\`u $D=[K(P):K]$, o\`u $g_0$ est la dimension du plus petit sous-groupe alg\'ebrique de $A$ contenant $P$ et o\`u $\kappa(g_0)=\left(2g_0(g_0+1)!\right)^{g_0+2}$.
\end{theo}
\demo C'est une cons\'equence imm\'ediate du corollaire 2 de \cite{ratazzi} appliqu\'e \`a la vari\'et\'e $V=\overline{\{P\}}$ image sch\'ematique de $P$ dans $A$ sur $K$. On peut faire une preuve directe (ce qui permet d'utiliser le r\'esultat principal de \cite{davidhindry} sans avoir \`a faire intervenir en plus leur remarque utilisant l'indice d'obstruction\footnote{remarque maintenant justifi\'ee par le pr\'esent article.}) : on commence par le cas o\`u $A=\prod_{i=1}^nA_i^{r_i}$, les $A_i$ \'etant des vari\'et\'es ab\'eliennes simples deux \`a deux non-isog\`enes et o\`u $\mathcal{L}$ est le fibr\'e en droites ample et sym\'etrique associ\'e au plongement 
\[A=\prod_{i=1}^nA_i^{r_i}\hookrightarrow \prod_{i=1}^n\mathbb{P}_{n_i}^{r_i}\overset{\textnormal{Segre}}{\hookrightarrow} \mathbb{P}_N,\]
\noindent les $A_i$ \'etant plong\'ees dans $\mathbb{P}_{n_i}$ par des fibr\'es en droites $\mathcal{L}_i$ tr\`es amples et sym\'etriques. On note $G$ le plus petit sous-groupe alg\'ebrique contenant $V$. On note $G^0$ la composante connexe de l'identit\'e de $G$. C'est une sous-vari\'et\'e ab\'elienne de $A$ et elle est donc isog\`ene \`a $B=\prod_{i=1}^nA_i^{s_i}$ o\`u $0\leq s_i\leq r_i$. On note alors $\pi : A \rightarrow B$ une projection naturelle obtenue par oubli de certaines coordonn\'ees, de sorte que $\pi_{\mid G}$ est une isog\'enie. Montrons que l'on est dans les conditions d'application du th\'eor\`eme principal de \cite{davidhindry} en prenant comme vari\'et\'e ab\'elienne $B$ et comme point $\pi(P)$.

\medskip

Si $\pi(P)$ est d'ordre fini modulo une sous-vari\'et\'e ab\'elienne stricte de $B$, en notant $H$ le plus petit sous-groupe alg\'ebrique contenant $\pi(P)$, on a $\textnormal{dim} H<\textnormal{dim} B$. Ainsi $G_1=G\cap\pi^{-1}(H)$ est un sous-groupe alg\'ebrique strict de $G$ (car $\pi_{\mid G}$ est une isog\'enie), contenant $V$. Ceci est absurde.

\medskip

Si $\pi(P)$ est d'ordre fini, comme $\pi$ est une isog\'enie, le point $P$ est aussi d'ordre fini. Ceci est absurde.

\medskip

\noindent Finalement, $\pi(P)$ est un point d'ordre infini modulo toute sous-vari\'et\'e ab\'elienne de $B$. On peut donc appliquer le th\'eor\`eme principal de \cite{davidhindry}. Par ailleurs, la hauteur et le degr\'e sont d\'efinis relativement aux plongements 
\[A=\prod_{i=1}^nA_i^{r_i}\hookrightarrow \prod_{i=1}^n\mathbb{P}_{n_i}^{r_i}\overset{\textnormal{Segre}}{\hookrightarrow} \mathbb{P}_{N}\ \ \textnormal{  et }\ \ B=\prod_{i=1}^nA_i^{s_i}\hookrightarrow \prod_{i=1}^n\mathbb{P}_{n_i}^{s_i}\overset{\textnormal{Segre}}{\hookrightarrow} \mathbb{P}_{N_B}.\]
\noindent De plus l'application $\overline{\pi} :  \prod_{i=1}^n\mathbb{P}_{n_i}^{r_i}\rightarrow  \prod_{i=1}^n\mathbb{P}_{n_i}^{s_i}$ est la projection lin\'eaire d\'efinie par oubli de coordonn\'ees. Dans ce cas, et pour ces plongements (en notant $\widehat{h}_{N_B}$ la hauteur de N\'eron-Tate de $B$ et $\widehat{h}_N$ celle de $A$ dans ces plongements), on a 
\[\widehat{h}_{N_B}(\pi(P))\leq \widehat{h}_N(P)\ \ \textnormal{ et }\ \ \deg\pi(P)\leq\deg P.\]
\noindent Ceci nous donne
\begin{align*}
\widehat{h}_N(P) 	& \geq \widehat{h}_{N_B}(\pi(P)),\ \ \ \textnormal{d'o\`u par le th\'eor\`eme de \cite{davidhindry},}\\
			& \geq \frac{c(B, N_B)}{\left(\deg\pi(P)\right)^{\frac{1}{g_0}}}\left(\log 2\deg\pi(P)\right)^{-\kappa(g_0)} \geq \frac{c(B, N_B)}{\left(\deg P\right)^{\frac{1}{g_0}}}\left(\log 2\deg P\right)^{-\kappa(g_0)}.\\
			& \geq \frac{c'(A, N)}{\left(\deg P\right)^{\frac{1}{g_0}}}\left(\log 2\deg P\right)^{-\kappa(g_0)},
\end{align*}

\medskip

\noindent o\`u on a pris pour $c'(A,N)$ le minimum des $c(B,N_B)$ quand $s_i$ varie dans $[\![0,r_i]\!]$.

\medskip

\noindent Dans le cas g\'en\'eral, la vari\'et\'e ab\'elienne $A$ est donn\'ee avec une isog\'enie $\rho$ vers la vari\'et\'e ab\'elienne $B=\prod_{i=1}^n A_i^{r_i}$. Soit $P$ d'ordre infini de la vari\'et\'e ab\'elienne de $A$. Le point $Q=\rho(P)$ est un point d'ordre infini de la vari\'et\'e ab\'elienne de $B$. Il r\'esulte facilement de la preuve de la proposition 14. de \cite{phi3} qu'il existe $c'(A,\mathcal{L})$ tel que
\[\widehat{h}_{\mathcal{L}}(P)\geq c'(A,\mathcal{L}) \widehat{h}_M(Q).\]
\noindent Ainsi en appliquant le r\'esultat pr\'ec\'edent, on en d\'eduit presque l'in\'egalit\'e voulue : il faut encore remplacer le degr\'e $\deg Q$ par $\deg P$. Or $\deg Q\leq \deg P$. Ceci permet de conclure.\hfill$\Box$

\medskip

\rem Ce r\'esultat am\'eliore le meilleur r\'esultat pr\'ec\'edemment connu, d\^u \`a Masser qui obtient dans \cite{lettre}, pour tout point $P$ d'ordre infini de $A(\overline{K})$ : 
\[\widehat{h}_{\mathcal{L}}(P)\geq \frac{c(A/K,\mathcal{L})}{D^2\log 2D}.\]

\medskip

\noindent En faisant la m\^eme preuve et en appliquant la partie (\ref{e1}) de la conjecture \ref{conj1app2} au lieu du th\'eor\`eme de \cite{davidhindry}, on obtient le

\medskip

\begin{cor}\label{co2} La partie (\ref{e1}) de la conjecture \ref{conj1app2} entra\^ine sa partie (\ref{e2}).
\end{cor}

\medskip

\rem Si au lieu du th\'eor\`eme de David-Hindry, on applique la conjecture \ref{conj2} sur le probl\`eme de Lehmer relatif, on peut partout remplacer le symbole $D$ par $D_{\tors}$ dans ce qui pr\'ec\`ede.

\subsection{Sur la conjecture de Lehmer mul\-tiho\-mo\-g\`ene sur les vari\'et\'es ab\'e\-lien\-nes}\label{multi}

\noindent Soit $A/K$ une vari\'et\'e ab\'elienne de dimension $g$. Quitte \`a augmenter un peu $K$ (\textit{cf.} par exemple \cite{ratazzi} lemme 1), on peut supposer (et on suppose) que tous les endomorphismes de $A$ sont d\'efinis sur $K$. On note $\widehat{h}_{\mathcal{L}}$ la hauteur de N\'eron-Tate sur $A(\overline{K})$ associ\'ee \`a un diviseur ample et sym\'etrique $\mathcal{L}$. Pour tout entier positif $n$ on pose $\mathcal{L}_n=\mathcal{L}^{\boxtimes n}$ fibr\'e en droites sym\'etrique ample sur $A^n$ et on note $\widehat{h}_{\mathcal{L}_n}$ la hauteur de N\'eron-Tate associ\'ee. On commence par un lemme.

\medskip

\begin{lemme}Soit $(P_1,\ldots,P_n)$ un point de $A^n(\overline{K})$. On a 
\[\widehat{h}_{\mathcal{L}_n}(P_1,\ldots,P_n)=\sum_{i=1}^n\widehat{h}_{\mathcal{L}}(P_i).\]
\end{lemme}
\demo C'est une cons\'equence formelle des propri\'et\'es de fonctorialit\'e des hauteurs de Weil et de la d\'efinition de la hauteur de N\'eron-Tate.\hfill $\Box$

\medskip

\noindent En utilisant ce lemme, on d\'emontre le r\'esultat suivant :

\medskip

\begin{theo}\label{theoapp2}Si $A/K$ est de type C.M., alors, pour tout entier $n\in \mathbb{N}$ il existe une constante $c(A/K,\mathcal{L},n)>0$ telle que pour tout point $(P_1,\ldots,P_n)\in A^n(\overline{K})$ d'ordre infini modulo toute sous-vari\'et\'e ab\'elienne stricte de $A^n$, on a :
\[\prod_{i=1}^n\widehat{h}_{\mathcal{L}}(P_i)\geq\frac{c(A/K,\mathcal{L},n)}{D^{\frac{1}{g}}}\left(\log 2D\right)^{-n\kappa(g)},\]
\noindent o\`u $D=[K(P_1,\ldots,P_n):K]$.
\end{theo}
\demo Soient $a_1,\ldots,a_n$ des entiers strictement positifs et $Q_1, \ldots, Q_n$ des points de $A(\overline{K})$ tels que pour tout $i$, $P_i=a_iQ_i$. On a 
\[\widehat{h}_{\mathcal{L}_n}(Q_1,\ldots,Q_n)=\sum_{i=1}^n\widehat{h}_{\mathcal{L}}(Q_i)=\sum_{i=1}^na_i^{-2}\widehat{h}_{\mathcal{L}}(P_i),\]
\noindent et,
\[\left[K(Q_1,\ldots,Q_n):K\right]^{\frac{1}{ng}}\leq\left(a_1^{2g}\times\cdots\times a_n^{2g}D\right)^{\frac{1}{ng}}.\]
\noindent Le th\'eor\`eme de David-Hindry nous donne alors
\[\sum_{i=1}^na_i^{-2}\widehat{h}_{\mathcal{L}}(P_i)\geq \frac{c(A/K,\mathcal{L},n)}{\left(\prod_{i=1}^na_i^{\frac{2}{n}}\right)D^{\frac{1}{gn}}}\left(\log \left((\prod_{i=1}^na_i)D\right)\right)^{-\kappa(g)}.\]
\noindent On pose maintenant, pour tout $1\leq i\leq n$,
\[ x_i=\frac{13\widehat{h}(P_i)}{4\min_j\widehat{h}(P_j)}, \text{ et } a_i=[\sqrt{x_i}].\]
\noindent Pour tout $i$, on a $x_i\geq \frac{13}{4}$ et $x_i\geq a_i^{2}\geq \frac{x_i}{3}$. Ainsi,
\[\sum_{i=1}^n a_i^{-2}\widehat{h}_{\mathcal{L}}(P_i)\leq \frac{3\times 4}{13}n\min_j\widehat{h}_{\mathcal{L}}(P_j),\ \ \text{ et }\ \ \prod_{i=1}^na_i^{2}\leq\left(\frac{13}{4\min_j\widehat{h}_{\mathcal{L}}(P_j)}\right)^n\prod_{i=1}^n\widehat{h}_{\mathcal{L}}(P_i).\]
\noindent Donc,
\begin{align*}
\min_j\widehat{h}_{\mathcal{L}}(P_j)	& \geq c_{10}(A/K,\mathcal{L},n)\sum_{i=1}^n a_i^{-2}\widehat{h}_{\mathcal{L}}(P_i)\\
				& \geq \frac{c_{11}(A/K,\mathcal{L},n)}{\left(\prod_{i=1}^na_i^{\frac{2}{n}}\right)D^{\frac{1}{gn}}}\left(\log 2D\prod_{i=1}^na_i\right)^{-\kappa(g)}\\
				& \geq \frac{4c_{11}(A/K,\mathcal{L},n)\min_j\widehat{h}_{\mathcal{L}}(P_j)}{13\prod_{i=1}^n\widehat{h}_{\mathcal{L}}(P_i)^{\frac{1}{n}}D^{\frac{1}{gn}}}\left(\log 2D\prod_{i=1}^na_i\right)^{-\kappa(g)}.
\end{align*}
\noindent Par ailleurs, on a la majoration
\[\log \prod_{i=1}^na_i\leq n\log\left(\frac{13}{2\min_j\widehat{h}_{\mathcal{L}}(P_j)}\right)+2\log \prod_{i=1}^n\widehat{h}_{\mathcal{L}}(P_i).\]
\noindent Or on peut toujours supposer que les $\widehat{h}_{\mathcal{L}}(P_i)$ sont inf\'erieurs \`a $1$, donc,
\[\log \prod_{i=1}^na_i\leq n\log\left(\frac{13}{2\min_j\widehat{h}_{\mathcal{L}}(P_j)}\right).\]
\noindent Ainsi, 
\[\log \left(2D\prod_{i=1}^na_i\right)\leq n\log\left(\frac{13D^{\frac{1}{n}}}{2\min_j\widehat{h}_{\mathcal{L}}(P_j)}\right).\]
\noindent On en d\'eduit que 
\[\prod_{i=1}^n\widehat{h}_{\mathcal{L}}(P_i)^{\frac{1}{n}}\geq \frac{c_1(A/K,n)}{D^{\frac{1}{ng}}}\left(\log\frac{D^{\frac{1}{n}}}{\min_j\widehat{h}_{\mathcal{L}}(P_j)}\right)^{-\kappa(g)}.\]
\noindent Le point $(P_1,\ldots,P_n)$ \'etant d'ordre infini modulo toute sous-vari\'et\'e ab\'elien\-ne, les points $P_i$ sont en particulier d'ordre infini sur $A$. Le r\'esultat inconditionnel de Masser sur la minoration de la hauteur des points sur les vari\'et\'es ab\'eliennes, theorem de \cite{masser}, nous donne donc :
\[\log \frac{D^{\frac{1}{n}}}{\min_j\widehat{h}_{\mathcal{L}}(P_j)}\leq c_2(A/K,\mathcal{L},n)\log 2D.\]
\noindent Ainsi, on en d\'eduit
\[\prod_{i=1}^n\widehat{h}_{\mathcal{L}}(P_i)\geq \frac{c_3(A/K,\mathcal{L},n)}{D^{\frac{1}{g}}}\left(\log 2D\right)^{-n\kappa(g)}\]
\noindent ce qui conclut.\hfill$\Box$

\medskip

\rem Si au lieu de faire appel au th\'eor\`eme 1.5. de \cite{davidhindry} dans la preuve du th\'eor\`eme \ref{theoapp2} on applique la conjecture \ref{conj1app2}, alors on en d\'eduit le r\'esultat suivant :

\medskip

\begin{theo}\label{theo2}Soient $A/K$ une vari\'et\'e ab\'elienne de dimension $g$ sur le corps de nombres $K$ et $\mathcal{L}$ un fibr\'e en droites sym\'etrique ample sur $A$. Si la conjecture \ref{conj1app2} est vraie pour $(A/K,\mathcal{L})$ alors, pour tout entier $n\in \mathbb{N}$ il existe une constante $c(A/K,\mathcal{L},n)>0$ telle que pour tout point $(P_1,\ldots,P_n)\in A^n(\overline{K})$ d'ordre infini modulo toute sous-vari\'et\'e ab\'elienne stricte de $A^n$, on a :
\[\prod_{i=1}^n\widehat{h}_{\mathcal{L}}(P_i)\geq\frac{c(A/K,\mathcal{L},n)}{D^{\frac{1}{g}}},\] 
\noindent o\`u $D=[K(P_1,\ldots,P_n):K]$.
\end{theo}

\medskip

\rem Les m\^emes remarques qu'au paragraphe A.1 pr\'ec\'edent, concernant le remplacement de $D$ par $D_{\tors}$ s'appliquent.

\medskip

\rem En fait dans leur article \cite{davidhindry}, les auteurs formulent \'egale\-ment une conjecture multihomog\`ene du probl\`eme de Lehmer ab\'elien. Plut\^ot que de supposer le point $(P_1,\ldots, P_n)$ d'ordre infini modulo toute sous-vari\'et\'e ab\'elienne stricte de $A^n$, ils supposent les points $P_i$ lin\'e\-aire\-ment ind\'ependants dans $A$. Pr\'ecis\'ement ils donnent la conjecture 1.6 suivante :

\medskip

\begin{conj}\label{conj2app2}\textnormal{\textbf{(David-Hindry)}} Soient $A/K$ une vari\'et\'e ab\'elienne de dimension $g$ sur un corps de nombres et $\mathcal{L}$ un fibr\'e en droites sym\'etrique ample sur $A$. Pour tout entier $n\in \mathbb{N}$ il existe une constante $c(A/K,\mathcal{L},n)>0$ telle que pour tout $n$-uplet $(P_1,\ldots,P_n)$ de points d'ordre infini dans $A(\overline{K})$, End($A$)-lin\'eairement ind\'ependants, on a :
\[\prod_{i=1}^n\widehat{h}_{\mathcal{L}}(P_i)\geq\frac{c(A/K,\mathcal{L},n)}{D^{\frac{1}{g}}},\] 
\noindent o\`u $D=[K(P_1,\ldots,P_n):K]$.
\end{conj}

\medskip

\noindent Dans la formulation de la conjecture \ref{conj2app2} qu'ils donnent, David-Hindry \'ecri\-vent ``li\-n\'e\-aire\-ment ind\'ependants'' sans pr\'eciser s'il s'agit de $\mathbb{Z}$-lin\'eairement ou de End($A$)-lin\'eairement ind\'ependants. Il para\^it pr\'ef\'erable de pr\'eciser. En effet, si on comprend l'assertion ``li\-n\'eaire\-ment ind\'ependants'' comme $\mathbb{Z}$-lin\'eaire\-ment ind\'e\-pen\-dants, alors la conjecture \ref{conj2app2} est fausse comme le montre l'exemple suivant : on prend $E/K$ une courbe elliptique \`a multiplication complexe par un corps quadratique imaginaire contenu dans $K$. On se donne $\alpha\in \textnormal{End}(E)$ un endomorphisme qui n'est pas la multiplication par un entier, on se donne \'egalement un point $P_1$ d'ordre infini dans $E(\overline{K})$ et pour tout $n\geq 1$, on choisit des points $P_n$ tels que $nP_n=P_1$. Enfin on pose $Q_n=\alpha(P_n)$. Puisque $P_1$ est d'ordre infini, les points $P_n$ et $Q_n$ sont $\mathbb{Z}$-lin\'eairement ind\'ependants. De plus on a 
\[\widehat{h}(P_n)\widehat{h}(Q_n)=\frac{\textnormal{N}(\alpha)}{n^4}\widehat{h}(P_1)^2,\ \ \text{ et }\ \ D_n:=[K(P_n,Q_n):K]=[K(P_n):K]\leq cn^2.\]
\noindent Donc,
\[\widehat{h}(P_n)\widehat{h}(Q_n)\leq \frac{c'}{D_n^2}.\]
\noindent Ceci montre que l'hypoth\`ese ``$\mathbb{Z}$-lin\'eairement ind\'ependants'' est insuffisante.

\medskip

\noindent Par contre en supposant les points  End($A$)-lin\'eairement ind\'ependants, la situation est bien meilleu\-re. Pr\'ecis\'ement, on a le

\medskip

\begin{theo}\label{theo3}La conjecture \ref{conj1app2} entra\^ine la conjecture \ref{conj2app2}.
\end{theo}
\demo Soit $n>0$ un entier. Au vu du th\'eor\`eme \ref{theo2}, la seule chose \`a prouver, est de montrer que l'hypoth\`ese (i) : ``les points $(P_1,\ldots,P_n)$ sont End($A$)-lin\'eairement ind\'ependants'', entra\^ine l'hypoth\`ese (ii) : ``le point $\mathbf{P}=(P_1,\ldots,P_n)$ est d'ordre infini modulo toute sous-vari\'et\'e ab\'elienne stricte de $A^n$.'' On va plut\^ot montrer que non(ii) implique non(i). Si non(ii) est vraie, alors, il existe un endomorphisme $\varphi$, non-nul, de $A^n$ tel que $\varphi(\mathbf{P})=0$. Or on peut \'ecrire $\varphi(\mathbf{P})=\left(\varphi_1(\mathbf{P}),\ldots,\varphi_n(\mathbf{P})\right)$, o\`u les $\varphi_i$ sont des morphismes de $A^n$ vers $A$ non tous nuls. On suppose par exemple que $\varphi_1$ est non-nul. En notant $\psi_i$ la restriction de $\varphi_1$ \`a la $i$-eme composante de $A^n$, on obtient ainsi $n$ endomorphismes de $A$, $\psi_1,\ldots,\psi_n$, non tous nuls et tels que 
\[\sum_{i=1}^n\psi_i(P_i)=\varphi_1(\mathbf{P})=0.\] 
\noindent Autrement dit, les points $P_1,\ldots,P_n$ sont End($A$)-lin\'eairement d\'ependants.\hfill$\Box$

\medskip

\noindent Enfin la m\^eme preuve permet de constater que le th\'eor\`eme \ref{theo2} entra\^ine un \'enonc\'e analogue en rempla\c{c}ant l'hypoth\`ese ``d'ordre infini modulo toute sous-vari\'et\'e ab\'elienne stricte'' par ``End($A$)-lin\'eairement ind\'ependants''. Ce dernier r\'esultat \`a \'egalement \'et\'e montr\'e par Viada \cite{viada} proposition 4. dans le cas particulier o\`u $A$ est une courbe elliptique.

\bigskip

\noindent \textbf{Adresse : }Nicolas Ratazzi\\
Universit\'e Paris-Sud 11\\ 
Batiment 425, Math\'ematiques\\
 91405, Orsay Cedex\\
 France

\end{document}